\documentstyle[amstex,amscd,12pt]{article}
\makeatletter
\date{}

\newtheorem{theorem}{Theorem}[section]

\newtheorem{proposition}[theorem]{Proposition}

\newcommand{\st}{{\rm st}}
\newcommand{\Int}{{\rm Int}}

\newcommand{\no}{$\rm  {N\ddot{o}beling}$ }
\newcommand{\noo}{$ {N\ddot{o}beling}$ }

\newcommand{\lo}{\longrightarrow}

\begin{document}

\title{A $Z$-set unknotting theorem for N\"{o}beling spaces\\
}

\author{Michael  Levin}

\maketitle
\begin{abstract}
We prove a   $Z$-set unknotting theorem for \no spaces. This
generalizes a result obtained by  S. Ageev   for a restricted
class of $Z$-sets. The theorem is proved for a certain model of
\no spaces.
\bigskip
\\
{\bf Math. Subj. Class.:} 55M10, 54F45.
\end{abstract}
\begin{section}{Introduction}\label{sec1}

 All spaces are assumed to be separable metrizable.
   A manifold means a manifold with (possibly empty) boundary and
   a triangulated space   means a locally finite simplicial
  complex which we identify with the underlying space.
  For a triangulated space we consider only triangulations
  compatible with the PL-structure of the space.
  All triangulated manifolds are assumed to be combinatorial.
 \\

 A complete $n$-dimensional metric space $X$ is  said to be an
 $n$-dimensional \no space if
 the following properties are satisfied:

 (i) $X $ is an absolute extensor in dimension $n$,
  that is every map
 $f: A \lo X$ from a closed subset $A$ of a space $Y$ of
 $\dim \leq n$
  extends over $Y$;

 (ii) every map $f : Y \lo X$ from a complete metric space
 $Y$ of $\dim \leq n$ can be arbitrarily closely
 approximated by a closed embedding, that is  for every open
 cover  $\cal U$ of $X$ there is a closed embedding $g : Y \lo X$
 which is $\cal U$-close to $f$
 ($\cal U$-close means that  for every $y \in Y$ there is
 an element of $\cal U$ that contains  both $f(y)$ and $g(y)$).
 \\

  Examples  of   \no spaces can be constructed as follows.

  By a rational map between triangulated spaces we mean a PL-map
  that sends points with rational barycentric coordinates
  to points with rational barycentric coordinates.
   Two triangulations  of a space   are said to be rationally
  equivalent if the identity map is a rational map with
  respect to these triangulations (it is easy to see that
  if a PL-homeomorphism is rational in one direction then it is
  rational in the opposite direction as well).
  Let $M$ be a triangulated space. Every triangulation of $M$
  which is rationally equivalent to the given triangulation
  of $M$ is said to be a rational triangulation and
  the class of all rational triangulations is said to be
  the rational structure of $M$.
   Denote by
  $M(k)$ the subspace of $M$ which is the complement of
  the union of all the triangulated spaces of $\dim \leq k$
  which are rationally  embedded in $M$.

  The space $M(k)$  admits the following
   interpretation. Fix a rational triangulation of $M$ and
  embed $M$ in the Hilbert space $H$  by an embedding  which sends
  the vertices  to points with rational coordinates and which is
  linear on every simplex of  the triangulation.  A $k$-dimensional plane
   in $H$ is said to be  rational if it is spanned by
   $k+1$  points with rational coordinates.
  Denote by $K$ the union  of all rational $k$-dimensional planes.
  Then $M(k)=M \setminus K$. Indeed, for every simplex $\Delta$
  of $M$ and every rational $k$-dimensional plane $L$,
  $\Delta  \cap L$ admits a triangulation for which it is rationally embedded
   in $M$ and therefore
  $M(k) \subset M \setminus K$. Now let $e: \Delta' \lo M$ be
  a rational embedding of  a simplex $\Delta'$ of $\dim \leq k$.
   Then there is  a (not necessarily rational)
    triangulation $\cal T$ of $\Delta'$ such that $e$ is linear
   on
   every simplex of $\cal T$. Since  every $\Delta'' \in \cal T$
   with $\dim \Delta''= \dim \Delta'$ has a dense subset of points
   with rational barycentric coordinates with respect to $\Delta'$
   we conclude that $e(\Delta'')$ is contained in a $k$-dimensional
   rational plane in $H$. Thus $e(\Delta') \subset K$ and hence
   $M(k)=M\setminus K$.

  Let us state the following important fact  leaving its proof to the
  reader.

 \begin{theorem}
 \label{t1} Let $M$ be a triangulated $m$-dimensional manifold,
  let $k\geq 0$ be an integer and
 let $n=m-k-1$.  If $M$ is $(n-1)$-connected and
 $m\geq 2n+1$ then $M(k)$ is an $n$-dimensional \noo space.

 \end{theorem}

  A space  $M(k)$ satisfying  the assumptions of Theorem \ref{t2}
 will be called a  \no space modeled on a triangulated manifold.

 A  subset $A$  of a space  $X$ is called a $Z$-set if
 $A$ is closed in $X$ and
 the identity map of $X$ can be arbitrarily closely approximated
 by a map $f : X \lo X$ with $f(X) \cap A =\emptyset$.
 Note that if $X$ is an $n$-dimensional \no space modeled
 on a manifold $M$
  and $A \subset X$ is a $Z$-set in $X$ then
 $X \setminus A$ is also an $n$-dimensional
  \no space modeled on the manifold $N=M \setminus$ the closure
  of $A$ in $M$ (the rational structure of $N$ is defined such that
  the inclusion is a rational map).

 The main result of this paper is
 \begin{theorem}
 \label{t2}
 Let $X_1$ and $X_2$ be $n$-dimensional \noo spaces and let $A_1$
 and $A_2$ be $Z$-sets in $X_1$ and $X_2$ respectively
 such that $X_1\setminus A_1$ and $X_2 \setminus A_2$  are
 homeomorphic to  $n$-dimensional \noo spaces modeled
 on triangulated manifolds.
 If $A_1$ and $A_2$ are homeomorphic
 then any homeomorphism  between $ A_1$ and $ A_2$ can be extended to
 a homeomorphism between $X_1$ and $X_2$.
 \end{theorem}

  Theorem \ref{t2} was proved  by Ageev \cite{a1} under
  some additional restrictions on the sets $A_1$ and $A_2$.
  The proof of Theorem \ref{t2} is self-contained and relies
  only on  well-known facts of PL-topology \cite{h,r}, \no spaces
  \cite{c} and elementary properties of partitions presented
  in the very beginning of \cite{b1}. Some ideas of
     the proof of Theorem \ref{t2} came from
     \cite{klt}.

    The results of    this paper along  with
     some additional arguments can be applied to
     validate  the characterization
     theorem  for \no
     saying  that any two \no spaces of the same dimension
     are homeomorphic. A proof of the characterization theorem
     will be presented in a forthcoming paper.
     Note that the characterization
  theorem implies that the assumptions on the complements
  $X_1\setminus A_1$ and $ X_2 \setminus A_2$ in Theorem \ref{t2} are
  automatically  satisfied and therefore can be be dropped.

     Let us finally mention that both a $Z$-set unknotting  theorem  and
     the characterization theorem for \no spaces were also
     announced by A. Nagorko at International Conference
     on Geometric Topology, July 2005, Bedlewo, Poland.

 \end{section}

 \begin{section}{Preliminaries}
 \subsection{General}\label{general}

 Let $M$ be a triangulated manifold.
 The geometric interior $\Int M$ of $M$  is
 the set of points  having  a neighborhood
 PL-homeomorphic to an Euclidean space of dimension=$\dim M$.
 The set $\partial  M =M \setminus \Int M$ is  the geometric boundary
 of $M$.

 A triangulated space (manifold)  PL-embedded in $M$
  is said to be a PL-subspace (PL-submanifold) of $M$.
 An open subset of  a triangulated space is always  considered
 with the induced PL-structure for which the inclusion is
 a PL-map.
  Note that an open subset of a
 PL-sudmanifold  of $M$  is also a PL-submanifold
  of $M$.

 A subset $P$ of $M$ is said to be a  PL-subcomplex
 of $M$ if there is a  triangulation  of $M$
 for which $P$ is a subcomplex. A subset $R$
 of $M$ is said to be a  PL-presented subset of $M$
  if there are closed subsets
  $R_1 \subset \dots \subset R_n$ of $M$ such that
  $R_n =R$,
 $R_1$ is a PL-subcomplex  of $M$
 and $R_{i+1} \setminus R_i$ is a PL-subcomplex
 of $M\setminus R_i$,  $i=1, \dots ,n-1$.

  A  collection $\cal P$ of subsets  of $M$ is said to be a decomposition
  of $M$ if
 $\cal P$ is a locally finite cover of  $M$
 and the elements of $\cal P$ are PL-subcomplexes of $M$.
  By a finite intersection of a decomposition
  ${\cal P}$ we   mean an intersection  of finitely many elements
  of ${\cal P}$ (the elements of ${\cal P}$ are
  also considered as finite intersections). Note that since ${\cal P}$ is locally finite,
  any non-empty intersection of elements of ${\cal P}$
  must be a finite intersection of $\cal P$.
  It is clear that  for a PL-submanifold $N$ of $M$,
  the restriction ${\cal P}| N=\{P \cap N : P \in {\cal P} \}$ of $\cal P$
  to $N$ is a decomposition of $N$.

   A decomposition  $\cal P$ of $M$ is said to be a partition of
   $M$ if
   each finite intersection   of $\cal P$ is
  a PL-manifold, the geometric interiors of any two non-equal
  finite intersections are disjoint and for every non-empty
  finite intersection $P=P_0\cap ...\cap P_t$ of
  distinct elements  $P_0,...,P_t  \in \cal P$,
   $\dim P =\dim M  -t$.

  A decomposition   $\cal P$   of $M$ is
  said to be a partition
  on
  a PL-submanifold $N$ of $M$ if  ${\cal P}| N$ is a partition of $N$
  (in this case we also say that $\cal P $ forms a partition on
  $N$ or $\cal P$ restricted to $N$ is a partition).

   A subset $P$ of $M$ is said to be $l$-co-connected, $l\geq 0$,
   if $P$ is $(\dim P -l)$-connected (we assume that $P$
   is $l$-co-connected for every $l \geq \dim P  +1$).
   A partition   $\cal P$ is said to be $l$-co-connected
    if
   every finite intersection  of $\cal  P $ is
   $l$-co-connected.

   A map  $f : X \lo M$  is said to be
 in general position  with a triangulation of $M$ if for every
 simplex $\Delta $ of the triangulation
 $\dim f^{-1}(\Delta) \leq \dim X +\dim \Delta -\dim M$.
 Every map from $X$ to $M$ can be arbitrarily closely approximated by
 a general position map. Moreover, let $F$ be closed in $X$
 and  a map $f : X \lo M$ restricted to $F$ be in general position.
 Then $f$ can be arbitrarily closely approximated
 by a general position map
 $g : X \lo M$ such that $g$ coincides
 with $f$ on $F$ and $g(X\setminus F) \subset \Int M $.
 A map $f : X \lo M$ is said to be in general position with a decomposition
 of $M$ if it is in general position with a triangulation
 underlying this decomposition (a triangulation of $M$ for which the elements
 of the decomposition are subcomplexes).

    By an $n$-dimensional
      cube
      $B^n$
      we mean a set of the form
      $B^n=\{ (x_1,...,x_n) : -r_i \leq x_i \leq r_i, r_i >0, i=1,...,n\}$
      in the Euclidean space $R^n$.
       Thus we always assume that the
      origin  $O$ is  at the center of  $B^n$. Considering
      the product $B^n \times B^m$ of two cubes we identify
      $B^n$ and $B^m$ with $B^n \times O$ and $O\times B^m$
      respectively.

      Let ${\cal A}$ and $  \cal B$ be collections of subsets of
      a set $X$ and let $C \subset X$. We denote
      $\st(C, {\cal A})=\cup \{ A : A \in {\cal A}, A \cap C \neq \emptyset \}$
      and $\st({\cal A}, {\cal B})=\{\st(A, {\cal B}): A \in {\cal A} \}$.

  \subsection{ Elementary properties of partitions}
   \label{facts}

   If  $\cal P$ is a  partition of a triangulated manifold $M$ then
   ${\cal P}|V$ is a partition of $V$ for every open $V \subset M$.
   Let $\cal P$ be  a decomposition of a triangulated manifold $M$ and let
   $\cal V$ be an open cover $M$ such that ${\cal P}|V$ is a partition
   of $V$ for every $V \in {\cal P}$ then
   $\cal P$ is a partition of $M$.

   Assume that $M$ is a triangulated space.
   One can show that   $M\times (0,1) $
   or $M \times [0,1)$  is
    a PL-manifold if and only if $M$
    is a PL-manifold.
    This implies that for a triangulated manifold $N$, the product
    $M\times N$ is a PL-manifold if and only if $M$ is a PL-manifold.
    Thus if $M$ and $N$ are triangulated manifolds and
      ${\cal P}$ is a decomposition of $M$
    then  the decomposition
    ${\cal P} \times N= \{ P \times N: P \in {\cal P}\}$
     is a partition of $M \times N$
     if and only if $\cal P$ is
    a partition of $M$.

  The  following  properties are proved in \cite{b1} for compact manifolds
  however their proof also applies for the non-compact case.

  Let $\cal P$ be a partition of a triangulated  manifold $M$,
  let ${\cal P}= \cup \{{\cal P}_i : i=1,2,\dots \}$ be a splitting
  of ${\cal P}$ into disjoint subfamilies ${\cal P}_i$
  and let $Q_i=\cup \{P : P\in {\cal P}_i \}$. Then
  ${\cal Q}=\{Q_1,Q_2, \dots\}$ is a partition of $M$ (see 1.1.5 of \cite{b1}).
  In particular, any union of elements of a partition is a PL-manifold.

    Let $\cal P$ be a partition of a triangulated manifold $M$.
    Then for every finite intersection $P$ of $\cal P$,
    $P \cap \partial M \subset \partial P$ (see 1.1.9 of  \cite{b1}).

    Let  $M$ be an $m$-dimensional triangulated  manifold, $\cal P$ a partition of $M$
    and $\cal Q$  a decomposition of $M$ such that for every finite intersection
   $P$  of $\cal P$, ${\cal Q}|P$ is a partition of $P$. Then
   $\cal Q$ is a partition of $M$ (see 1.1.11 of \cite{b1}). In particular,
   if $M=M_1\cup M_2$ is a decomposition of $M$ into two $m$-dimensional
   PL-submanifolds $M_1$ and $M_2$ such that $N=M_1\cap M_2$ is
   an $(m-1)$-dimensional PL-submanifold of both $\partial M_1$ and $\partial M_2$ and,
   $\cal Q$ is a decomposition of $M$ such that ${\cal Q}|M_1$,
     ${\cal Q}|M_2$ and ${\cal Q}|N$ are partitions then $\cal Q$ is
     a partition of $M$.

  \subsection{ A Matching of partitions}\label{match}

  Let ${\cal P}_1$ and ${\cal P}_2$ be   partitions
  of triangulated  manifolds $M_1$ and $M_2$ respectively.
  A one-to-one correspondence  $\mu : {\cal P}_1 \lo {\cal P}_2$ is said to be
  a matching  of ${\cal P}_1$ and ${\cal P}_2$ if
  it induces a one-to-one correspondence between non-empty
   intersections.
   This  means that for  every finite intersection
   $P =P_0\cap \dots \cap P_t$ of distinct elements
   $P_0,...,P_t  \in {\cal P}_1$,
    $\mu(P)=\mu(P_0)\cap \dots \cap\mu( P_t)=\emptyset$
    if and only if $P=\emptyset$.

  Assume that  $\mu : {\cal P}_1 \lo {\cal P}_2 $ is a matching of
  partitions ${\cal P}_1$ and ${\cal P}_2$ such that
  ${\cal P}_2$ is $l$-co-connected.
    Let  $F$ be a closed
    subset of $X$  such that $\dim X \leq m_i-l +1$,
   $m_i=\dim M_i$ $i=1,2$,
    and let   $f_i : F \lo M_i$ be maps such that
    $f^{-1}_1(P)=f^{-1}_2(\mu(P))$  for
    every finite intersection  $P$ of ${\cal P}_1$.
    Assume that a map $f'_1 : X \lo M_1$ extends $f_1$  such that
    $f'_1$ is in general position with ${\cal P}_1$.
    Then there is a map $f'_2 : X \lo M_2$ such that $f'_2$ extends
    $f_2$ and $(f'_1)^{-1}(P)=(f'_2)^{-1}(\mu(P))$  for
    every finite intersection  $P$ of ${\cal P}_1$.

     The required extension of $f_2$ can be constructed by
    induction on
    the co-dimension $t= m_2, m_2-1,m_2-2,\dots, 0$
    of the intersections of ${\cal P}_2$. Assume that we already
    extended $f_2$ to $f'_2$ on $X_t=$the union
    of the preimages of all  the intersections of ${\cal P}_2$
    of $\dim \leq m_2- t$  such that $f'_2$ on $X_t$
    satisfies the required properties
    and let an intersection $P$ of
    ${\cal P}_2$  be of $\dim = m_2 -t+1$. Then  for $X_P
    =(f'_1)^{-1}(\mu^{-1}(P))$,
    $\dim X_P \leq \dim X +\dim \mu^{-1}(P) -m_1=
    \dim X +\dim P - m_2$. Since $P$ is $(\dim P -l)$-connected
    and $\dim P - l  \geq \dim P -m_2+\dim X -1\geq
    \dim X_P -1$,  $f'_2$ restricted to $X_P \cap X_t$ can be
    extended over $X_P$ such that $f'_2(X_P \setminus X_t)\subset
    \Int P$.

    Thus we can extend $f_2$ to a map $f'_2$ with the required
    properties.  This extension
    will be called a transfer of the extension $f'_1$
    via the matching $\mu$.

    Assume that $\mu : {\cal P}_1 \lo {\cal P}_2$ is
    a matching of  partitions on manifolds  $M_1$ and
    $M_2$ respectively such that ${\cal P}_i$ is
    $l_i$-co-connected and $m_1-l_1 =m_2-l_2$ where
    $m_i =\dim M_i$. Let $n\leq m_1- l_1$ and
    let $P_1 \subset P'_1 $ be
    subsets of $M_1$ such that $P_1, P'_1$ are  unions of
    elements of ${\cal P}_1$
    and  the inclusion $P_1 \subset P'_1$ induces
    the zero-homomorphism of the homotopy groups in $\dim \leq n$.
    Then the inclusion $P_2=\mu(P_1) \subset P'_2 =\mu(P'_1)$
     also induces
    the zero-homomorphism of the homotopy groups in $\dim \leq n$.
    Indeed, take
    a map $f_2 : S^p \lo P_2$, $p \leq n$ from a $p$-dimensional sphere $S^p$
    into $P_2$. Since $P_2$ is a manifold
    we can homotope  $f_2$ inside $P_2$ to a general position map.
    Thus we assume that $f_2$ is a general position map
    and transfer this map via the matching $\mu$ to
    $f_1 : S^p \lo P_1$. Extend $f_1$ to a general position map
    $f'_1  : B^{p+1} \lo P'_1$ of an $(p+1)$-dimensional ball $B^{p+1}$
    such that $S^p =\partial B^{p+1}$ and once again
    transfer this extension to an extension
    $f'_2 : B^{p+1} \lo P'_2$ of $f_2$. Thus $f_2$ is
    null-homotopic in $P'_2$.

  \subsection{A black hole modification }\label{black}

  Let $ B^q, B^{m-q}, B^{m-q}_*$ and $B^m= B^q \times B^{m-q}$ be
  cubes with the dimensions indicated by the superscripts
   such that $B_*^{m-q}$ is contained in
   $\Int B^{m-q}$ and the centers of $B_*^{m-q}$ and $B^{m-q}$ coincide
  (recall that we assume that the centers of the cubes
  are located at $O$ and we identify $B^q$ and $B^{m-q}$  with the subsets
  $B^q\times O$ and $O \times B^{m-q}$ of $B^m$ respectively).
  Denote
  $T=B^q \times (B^{m-q}
  \setminus \Int  B^{m-q}_*)$ and
  $S= \partial B^{m-q}$.
  Note that $S$ is an $(m-q-1)$-dimensional sphere
  in $ \partial T$ and $T$ is a subset of $B^m$
  which is  PL-homeomorphic to
  the product  $B^{q+1} \times S$ where $S$ is identified with
  $O \times S$.

   Assume that $B^m$ is PL-embedded
   in  $\Int M$ of  a triangulated $m$-dimensional manifold $M$.
   Denote $L= M\setminus \Int T $ and assume
   that $F_L$ is a PL-subcomplex of $L$ and
   ${\cal P}_L$  is a  decomposition of $L$ such
   that $S \subset F_L$ and
   ${\cal P}_L$ restricted to $L\setminus F_L$ and
  to  $\partial T \setminus F_L$
    are partitions.

   Fix a triangulation $\cal T$ of $M$ that underlies $S$, $F_L$
   and the elements
   of ${\cal P}_L$ and denote by $\cal A$ the collection
   of subcomplexes of $L$ with respect to this triangulation.
   We are going to extend each   $A$ in $\cal A$  to
   a PL-subcomplex  $A\subset A^\pi$ of $M$  such that
   the following conditions are satisfied:

   ${\cal P}_L^\pi =\{ P_L^\pi : P_L \in {\cal P}_L \}$ is a decomposition of
   $M$ which forms a partition on
   $M \setminus   F_L^\pi$;

   $F_L^\pi \cap L =F_L$ and ${\cal P}^\pi_L$ restricted to $L \setminus F^\pi_L$
   coincides with ${\cal P}_L$ restricted to $L \setminus F_L$;

   the correspondence between ${\cal P}_L$ and ${\cal P}_L^\pi$
    defined by $P_L \lo P_L^\pi$ induces a matching
    of partitions when  ${\cal P}_L$ and ${\cal P}_L^\pi$ are
    restricted to $L\setminus F_L$ and $M\setminus F_L^\pi$
    respectively;

    for every intersection $P_L= P^0_L\cap \dots \cap P^j_L$
    of  elements
    $P^0_L,\dots, P^j_L$ of ${\cal P}_L$,
    $P_L^\pi \setminus F_L^\pi =
    ((P^0_L)^\pi \cap \dots \cap (P^j_L)^\pi)\setminus F_L^\pi$
    and $P_L \setminus F_L$ is a (strong) deformation retract
    of $P_L^\pi \setminus F_L^\pi$;

    $\dim (F_L^\pi \cap T)\setminus S \leq
    \dim (F_L\cap \partial T)\setminus S +1$;

    and finally $F_L^\pi =S$ if $F_L=S$.\\

   Note that $T \setminus S$ can be represented as the product
   $(\partial T \setminus S) \times [0,1)$ where $\partial T
   \setminus S$ corresponds to
   $(\partial T \setminus S)\times 0$.
   One of the ways to construct $A^\pi$ is to find  a retraction
      $\pi : T\setminus S \lo \partial T \setminus S$
    which is topologically equivalent
   to the projection of  $(\partial T \setminus S) \times [0,1)$
   to $(\partial T \setminus S)\times 0$ such that for  $A^\pi$ defined
   by $A^\pi =A$ if $A\cap \partial T=\emptyset$ and
   $A^\pi = A\cup S \cup \pi^{-1}(A\cap \partial T)$ otherwise we
   have that $A^\pi$
   is a  PL-subcomplex of $M$ for every $A \in \cal A$.
   Such a retraction  $\pi$ can be constructed however the author
   could not find an elementary and direct argument for this. Therefore
    for the sake of elementariness (and at expense of a geometric
    transparency) we adopt a slightly different approach for constructing $A^\pi$
     which is
    presented in \ref{ext}. Note that all  the  properties of ${\cal P}_L^\pi$
    and $F_L^\pi$ that will be used in the sequel are described
    above. Thus for understanding  the rest of the paper the
    reader  may assume that the construction is carried out on
    base of an appropriate retraction $\pi$.

  We will call  the decomposition ${\cal P}_L^\pi$  a black hole
   modification of ${\cal P}_L$ with the sphere $S$ as the
   black hole  of the modification.

 \subsection{ Improving connectivity of an intersection}\label{connect}
  Assume that $m \geq 2q +1$ and  $l=m-q+2$,  $M$ is
  a triangulated  $m$-dimensional
 $(q-1)$-connected (=$(l-1)$-co-connected) manifold and
  $F$ is a PL-subcomplex of $M$ lying in
   $\Int M $. Let
 $\cal P$ be a decomposition of $M$ which forms  an  $l$-co-connected
  partition on $U=M\setminus F$
  and let $P=P_0\cap P_1\cap \dots\cap P_t$, $0\leq t \leq m-l+1$
  be an  intersection of distinct elements
  $P_0,\dots, P_t$ of $\cal P$ such that $P\cap U\neq\emptyset$. Then
    $P\cap U$ is an $(m-t)$-dimensional and
  $(q-t-2)$-connected (= $l$-co-connected) manifold.

  Assume that
   the intersections of  ${\cal P}|U$  of
  dimension$>m-t$ are  $(l-1)$-co-connected. Let us show how
  to improve the connectivity of $P\cap U$ preserving the level of
   connectivity
  of the rest of the intersections.

  Let  $f : S_P \lo \Int (P\cap U) $
  be a PL-embedding of
  a triangulated $(q-t-1)$-dimensional sphere $S_P$
  such that  $f$
  is not null-homotopic in $P\cap U$.

   Consider a $t$-dimensional simplex $\Delta$ with
  vertices $V(\Delta)=\{v_0,..., v_t \}$.
  Let $\Delta'$  be a (sub)simplex of $\Delta$ lying in
   $\partial \Delta$ and let
  $V(\Delta') \subset V(\Delta)$ be the set of vertices
  of $\Delta'$.
    Denote $P(\Delta')=\cap \{P_i : v_i \in V(\Delta)\setminus
    V(\Delta')\}$.  Note that $P(\emptyset)=P$ and
    $P(\Delta') \subset P(\Delta'')$ if $\Delta' \subset
    \Delta''$.

   For every simplex $\Delta'$ of $\partial \Delta$ we
   will define by induction on $\dim \Delta'$
  a PL-embedding $f_{\Delta'} : S_P * \Delta' \lo  P(\Delta')\cap U$
  such that $f_{\Delta''}$ extends $f_{\Delta'}$ if $\Delta'
  \subset \Delta''$. Set $f_{\emptyset}=f$ and
  we are done for $t=0$.
  Let  $t >0$.
   Take a simplex $\Delta '$ of $\dim=i+1 \leq t-1$ and assume that
  $f_{\Delta'}$  is defined for all $\Delta'$ of $\dim \leq i$.
  Then the maps
  $f_{\Delta''}$ for $\Delta'' \subset \partial \Delta'$ define
  the corresponding PL-embedding $f_{\partial \Delta'} : S_P * \partial \Delta' \lo
  \partial (P(\Delta')\cap U)$. Note that $S_P * \partial \Delta'$ is a
  $(q-t+i)$-dimensional sphere and $P(\Delta')\cap U$ is
  $(m-t+i+2)$-dimensional and
  $(q-t+i+1)$-connected (=$(l-1)$-co-connected) and
  $\dim (P(\Delta')\cap U)=m-t+i+2\geq 2q+1-t+i+2
  \geq 2(q-t+i+1)+ 1 +(t-i-2) \geq 2(q-t+i+1) +1$
  since $m\geq 2q+1$ and $t-i-2\geq 0$. Then $f_{\partial \Delta'}$ can be
  extended to a PL-embedding $f_{\Delta'} :S_P * \Delta' \lo
  P(\Delta')\cap U$ such that $f_{\Delta'} $ sends  $\Int(S_P * \Delta')$
  to  $\Int(P(\Delta')\cap U)$.

  The maps $f_{\Delta'}$ define the PL-embedding
  $f_{\partial \Delta}$ of $ S_P * \partial \Delta $ to
   $ \Int((P_0 \cup ...\cup P_t)\cap U)$.
  Identify the $(q-1)$-dimensional sphere $S^{q-1}=S_P * \partial \Delta$
  with its  image under $f_{\partial \Delta}$ and,
  by \ref{hole1},
  extend  the embedding of $S^{q-1}$ to a PL-embedding
   of a  $q$-dimensional cube
  $B^q$  such that $S^{q-1}=\partial B^q$ and
  the embedding of $B^q$ extends to a PL-embedding of an $m$-dimensional cube
  $B^m =B^q \times B^{m-q}\subset  \Int M$ having the following properties:
  $P' \cap B^m=(B^q \cap P') \times B^{m-q}$ for every $P'\in {\cal P}$;
  $F \cap B^m =(B^q\cap F)\times B^{m-q}$;
  $\cal P$ restricted to $N \cap U$ and $\partial B^m \cap U$ are partitions
  where  $N=M\setminus \Int B^m$ and
   for every finite intersection $P'$
  of  $\cal P$ such that $P'\cap U\neq \emptyset$
  we have that $P'\cap N \cap U \neq\emptyset$ and
      $P'\cap N \cap U $  is
  a deformation retract of $(P'\cap U) \setminus \Int   B^q$.
   Note that
   for every
   finite intersection $P'$
   of ${\cal P}$ that intersects $U$,
   $(\dim(P'\cap U) -(l-1) + 1)+ q+1=
   \dim(P' \cap U) -m+ 2q +1 \leq \dim (P'\cap U)$
   and therefore
    the inclusion
    $(P'\cap U)\setminus  \Int B^q
  \subset P'\cap U$ induces an isomorphism of
  the homotopy groups in  $\dim \leq \dim (P'\cap U) - (l-1)$
  (= in co-dimensions $\geq l-1$).
  Also note that since $F \cap \partial B^q =\emptyset$ we have
  $F \cap \partial B^m =(F \cap \Int B^q)\times \partial B^{m-q}$
  and therefore $\dim F \cap \partial B^m \leq \dim F -1$.

  Let $B^{m-q}_*$ be a cube lying in the geometric interior of
  $B^{m-q}$ such that the centers of $B^{m-q}_*$ and $B^{m-q}$
  coincide and
  $\partial B^q  \times B^{m-q}_* \subset \Int((P_0\cup \dots\cup P_t)\cap U)$.

   Use the notations of \ref{black}. Define
      $F_L =(N\cap F)\cup S$ where $S= \partial B^{m-q}$.
  Consider $B^q$ as the join $O* \partial B^q$ where $O$ is
  the center of $B^q$ and
  define a decomposition ${\cal P}_L$ of $L$ as  the collection
  of the sets ${\cal P}_L=\{ P'_L : P' \in {\cal P}\}$ where
  $P'_L=P'\cap N$ if $P' \in \cal P$ and $P'\neq P_0,\dots,P_t$
  and
  $P'_L = (P'\cap N) \cup ((O *(P' \cap \partial B^q )) \times B^{m-q}_*)$
  for $P'=P_i$, $0 \leq i \leq t$. Note that  ${\cal P}_L$ restricted to
  $N$ coincides with ${\cal P}$ restricted to $N$. It is easy to see that
  ${\cal P}_L$ restricted to $B^q \times B_*^{m-q}$,
  $\partial B^q \times (B^{m-q}\setminus \Int B_*^{m-q})$,
  $\partial B^q \times B_*^{m-q}$,
  $B^q \times \partial B_*^{m-q}$ and
   $\partial B^q \times \partial B_*^{m-q}$
   are partitions. Then, by \ref{facts},
   ${\cal P}_L$ is a partition on $L \setminus (N\cap F)$ and
  $\partial T \setminus (N\cap F)$  and hence
   ${\cal P}_L$ is a partition on $L \setminus F_L$
   and $\partial T \setminus F_L$.
  For a finite intersection  $P'=P'_0 \cap \dots \cap P'_j$ of
   $P'_0, \dots , P'_j\in \cal P$ such that $P' \cap U \neq \emptyset$
  denote
  $P'_L=(P'_0)_L \cap \dots \cap (P'_j)_L$.

  Now apply the black hole modification to ${\cal P}_L$ and $F_L$
  with
  $S$ as  the black hole of the modification.
    Note that since $S$  lies in $\partial N$ then, by \ref{facts},
  for every finite intersection $P'$ of ${\cal P}$,
  $P'\cap U\cap N \cap S \subset \partial (P'\cap U\cap N)$
  and hence
  the inclusion $(P'\cap U\cap N) \setminus S \subset P'\cap U\cap N$
  induces an isomorphism of the homotopy groups in
  all dimensions.
  Recall that
  $f(S_P) \subset  P\cap U \cap N$.
  By \ref{black} and the construction of ${\cal P}_L$ one can
  easily verify that
  for every finite intersection
  $P'$ of ${\cal P}$ such that such that
  $P'\cap U\neq \emptyset$ and $P \neq P'$,
  $(P'\cap U \cap N)\setminus S$ is a deformation retract of
  $(P'_L)^\pi \setminus  F_L^\pi$ and for the intersection $P$
  the inclusion of $(P\cap U\cap N)\setminus S$
   into $P_L^\pi \setminus F_L^\pi$ induces an
  isomorphism of the homotopy groups in dimensions $< q-t-1$ and an
  epimorphism in $\dim =q-t-1$ such that the element of the homotopy group
  represented by  the map $f$
  is in the kernel of this epimorphism.

  Thus for every intersection $P'$ of ${\cal P}$ such that
  $P'\cap U \neq \emptyset$
  the connectivity of $(P'_L)^\pi \setminus F_L^\pi$
  fits the connectivity of $P'\setminus F$ for $P' \neq P$ and
  we contributed toward improving the connectivity of
  $P \setminus F$ in its modification
  $P_L^\pi \setminus F^\pi_L$.

  Denote $F^\pi =F^\pi_L$, $(P')^\pi =(P'_L)^\pi$ for
  a finite intersection
  $P'$ of  $ {\cal P}$ and let ${\cal P}^\pi =
  \{(P')^\pi : P' \in {\cal P}\}= {\cal P}_L^\pi$.
 It is easy to see that
   the  one-to-one correspondence
   between ${\cal P}$ and ${\cal P}_L$ defined
   by $P' \lo (P')_L, P' \in {\cal P}$ becomes a matching
   of partitions when ${\cal P}$ and ${\cal P}_L$ are restricted
   to $M \setminus F$ and $L\setminus F_L$
   respectively.
  Then by \ref{black}
   the  one-to-one correspondence
   between ${\cal P}$ and ${\cal P}^\pi$ defined
   by $P' \lo (P')^\pi, P' \in {\cal P}$ becomes a matching
   of partitions when ${\cal P}$ and ${\cal P}^\pi$ are restricted
   to $M \setminus F$ and $M\setminus F^\pi$
   respectively. Note that $F^\pi \subset \Int M $ and
   because
   $\dim (F\cap \partial B^m )\leq \dim F   -1$ we have
    $\dim F^\pi \leq \max \{\dim S, \dim F\}=
    \max \{ m-q-1, \dim F\}$.

    Finally note that
    even if  $F=\emptyset$
    we have $F^\pi=S$ and
    we improve the connectivity of $P$
   at expense of
   an "irregular" behavior  of ${\cal P}^\pi$ on the black
  hole $S$, in particular, we create additional intersections
  of ${\cal P}^\pi$ on $S$ and
   the elements of ${\cal P}^\pi$ are no
  longer manifolds in $M$.

 \subsection{ A discretization of maps' images} \label{desc}
 Assume that  $M$ be a triangulated manifold  such that
 $\dim M \geq 2r+1$ and  $\cal U$ is an open cover  of $M$.

 Let $\cal F$ be countable collection of maps
  $f:S^{r} \lo M $ from
 a triangulated $r$-dimensional sphere $S^{r}$.
 Then there is a
 PL-subcomplex $R$ of $M$ of $\dim \leq r$
 lying in $\Int M $ such that
 every  map $f$  in ${\cal F}$ admits
 a $\cal U$-close approximation
 $f' : S^{r} \lo \Int M \setminus R $
 such that $f'$ is a PL-embedding and
 the images of  the maps in
 ${\cal F}'=\{ f' : f\in {\cal F }\}$ form
 a discrete family in $M\setminus R$.
   Moreover, given a decomposition $\cal P$ of $M$, the set $R$ in
    can be chosen so that $R$  is nowhere dense on the finite
    intersections of
     $\cal P$, that is the intersection of $R$ with every finite intersection
     of $\cal P$ is nowhere dense in that intersection.

  Indeed,
  arrange $\cal F$ into a sequence
 ${\cal F}=\{ f_1,f_2,\dots \}$ and
 take a sufficiently small triangulation $\cal T$
 of $M$ such that $\cal T$ underlies $\cal P$.
 Approximate each  $f_i : S^{r-1} \lo M$
  by a map $f''_i : S^{r} \lo M$ such that
 $f''_i(S^{r})$ is contained in  $R'=$the $r$-skeleton of $\cal T$.
 Fix a metric $d$ on $M$ and approximate each $f''_i$ by an PL-embedding
 $f'''_i : S^r \lo  M$
 such that $f'''_i(S^r) \cap R =\emptyset$, $ d(f'''_i(S^r), R')\leq 1/i$ and
 $f'''_i (S^r)\cap f'''_j (S^r)= \emptyset$ for $i\neq j$. Then
 the sets $f'''_i(S^r),i=1,2\dots$ form a discrete family in
 $M \setminus R' $. Now take a PL-embedding
 $g : M \lo  M$  such that $g(M)$ is a PL-subcomplex of $M$,
 $g(M)\subset \Int M$, the restriction of $g$ to $R'$  is in general position with
 $\cal T$ and $g$ is  close to the identity map of $M$.
 It is clear  that the construction  can be carried out  such that
 $f'_i=g \circ f'''_i$ is $\cal U$-close to $f_i$. Then
 the maps $f'_i$ and the set $R=g(R')$ will have the required properties.

 In a similar way one can show that if $\cal F$ is
 a countable collection of maps
  $f:B^r \lo M $ from
 a triangulated $r$-dimensional ball $B^r$
 such that for every $f$ in ${\cal F}$, $f$ restricted to $\partial B^r$
 is a PL-embedding and the images of the maps in
 $\{ f|_{\partial B^r} : f \in {\cal F}\}$
 form
 a discrete family in $M$
 then there is
 a PL-subcomplex $R$ of $M$ of $\dim \leq r$
 lying in $\Int M $ such that
 every  map $f$  in ${\cal F}$ admits
 a $\cal U$-close approximation
 $f' : B^r \lo M \setminus R $
 such that $f'$ is a PL-embedding,
 $f'|_{\partial B^r}=f|_{\partial B^r}$,
 $ f'(\Int B^r ) \subset\Int  M $
  and  the images of  the maps in
 ${\cal F}'=\{ f' : f\in {\cal F}\}$ form
 a discrete family in $M\setminus R$.
   Moreover, given a decomposition $\cal P$ of $M$, the set $R$ in
    can be chosen so that $R$ is nowhere dense on the finite intersections
     of $\cal P$.

   Note that if $M$ is $(r-1)$-connected
   and   $R$ is
  a closed PL-presented  subset of $M$ of $\dim \leq r$
   then $M'=M\setminus R$ is $(r-1)$-connected.
   In addition, if $\cal P$ is a decomposition of $M$
   which is an $l$-co-connected partition on an open subset
   $U$ of $M$ with $l \geq r+2$ and $R$ is nowhere dense
   on the finite intersections of $\cal P$ then
   ${\cal P}'=\{ P \setminus R : P \in {\cal P}\}$
   is a decomposition of $M'$ which is an
   $l$-co-connected partition on $U' =U\setminus R$
   and
   the one-to-one correspondence between
   $\cal P$ and ${\cal P}'$ defined
   by $P \lo P \setminus R, P \in \cal P$ induces a matching
   of partitions when $\cal P$ and ${\cal P}'$
   are restricted to $U$ and $U'$ respectively
   (the assumption that $R$ is nowhere dense on the finite
   intersections of ${\cal P}$ is needed to guarantee that
   for every  finite intersection $P$ of $\cal P$,
    $(P\cap U) \setminus R \neq \emptyset$ provided
   $P\cap U \neq \emptyset$).

  \subsection{ Improving connectivity of
  intersections simultaneously  }\label{simul}
  Adopt the notations  and  the assumptions
  of \ref{connect}.  We are going to show how to improve simultaneously
  the connectivity  of all the intersections of ${\cal P}| U$ of $\dim=m-t$
  to $(l-1)$-co-connectivity.

  For an intersection
   $P$ of $\cal P$
   with  $\dim P\cap U  =m-t$ choose a countable collection ${\cal F}_P$
   of PL-embeddings
   $f :  S_P \lo \Int(P\cap U )$
   from
   a  $(q-t-1)$-dimensional sphere $S_P$ that generate
   the $(q-t-1)$-dimensional homotopy group of $P\cap U$.
   Since $F$ is a PL-subcomplex of $M$ we may assume
   that the images of the maps in ${\cal F}_P$ lie
   outside a neighborhood of $F$. Then by \ref{desc}
   one can find a PL-subcomplex $R_P$  of $M$  such that
   $R_P \subset \Int(P\cap U) $ and
    $\dim R_P \leq q-t-1$
   and assume that the images of the embeddings in ${\cal F}_P$
    are contained in $\Int ((P\cap U )\setminus R_P)$ and
    form a discrete family in $M\setminus R_P$.

   Recall that in \ref{connect} the elements of $\cal P$ that form
  the intersection $P$ are enumerated according to the vertices
  of a sample simplex $\Delta =\{ v_0,\dots,v_t \}$ (we fix such enumeration
  arbitrarily and independently for every intersection $P$ of $\cal P$
  with $\dim P \cap U =m-t$).
  Following \ref{connect}
  and using  \ref{desc} we can  define for every
  $\Delta' \subset   \Delta$
  a set  $R_P^{\Delta'}$ and  collections of maps
  ${\cal F}_P^{\Delta'}=\{ f_{\Delta'}: f \in {\cal F}_P \}$ and
   ${\cal F}_P^{\partial \Delta'}=\{ f_{\partial \Delta'}: f \in {\cal F}_P \}$
   such that $R^{\Delta''}_P \subset R^{\Delta'}_P$
  if $ \Delta'' \subset \Delta'$,
   $R^{\Delta'}_P$ is a closed subset of $M$ lying in $P(\Delta')\cap U$,
 the images of the
  maps of ${\cal F}_P^{\Delta'}$ are contained in
  $G_P^{\Delta'}= (P(\Delta')\cap U)\setminus R_P^{\Delta'}$ and form
  a discrete family in $M \setminus R_P^{\Delta'}$,  and $R_P^{\Delta'}
  \setminus R_P^{\partial \Delta'}$
  is a PL-subcomplex of
  $M\setminus R_P^{\partial \Delta'}$ lying in
   $(P(\Delta')\cap U)\setminus R_P^{\partial \Delta'}$
  such that  $ \dim R_P^{\Delta'} \setminus R_P^{\partial
  \Delta'}\leq q-t +\dim \Delta'$
  where $R_P^{\partial \Delta'}=\cup\{ R^{\Delta''}_P :
  \Delta''\subset \partial \Delta' \}$.
  Here we assume that  $\Delta'=\emptyset$ is
   a simplex of $\Delta$ and  ${\cal F}_P^\emptyset = {\cal F}_P$
   and $R_P^\emptyset =R_P$ and we also assume that $P(\Delta)=M$.

   Denote by ${\cal F}^{\partial \Delta}$ and $R^{\partial \Delta}$
   the union of ${\cal F }_P^{\partial \Delta}$ and
   $R_P^{\partial \Delta}$ respectively over all the
   intersections $P$ of $\cal P$ with $\dim P \cap U =m-t$.
   In order to carry out
   the  construction described above such that
   $R^{\partial \Delta}$ will be a closed PL-presented subset
   of $M$ and the images of the maps in
   ${\cal F}^{\partial \Delta}$
   will form a discrete family in $M \setminus R^{\partial \Delta}$,
    we need to take into account
      that the same intersection $ P'$ of
   $\cal P$  with  $\dim  P' \cap U > m-t$ can be involved in many (even
   countably many) intersections of $\dim=m-t$. It can be done as follows.
   For a finite intersection  $P'$  of  ${\cal P}$ with
    $\dim P' \cap U > m-t$
    denote by  $R_{ \partial P'}$, $R_{P'}$,
    ${\cal F}_{ \partial P'}$
    and ${\cal F}_{ P'}$
     the union of $R_P^{ \partial \Delta'}$, $R_P^{\Delta'}$,
    ${\cal F}_P^{\partial \Delta'}$ and
    ${\cal F}_P^{\Delta'}$ respectively over all finite intersections $P$ of $\cal P$
   such that $\dim P \cap U=m-t$ and $\Delta' \subset \partial \Delta$
   such that
    $P(\Delta')=P'$ and, let $G_{P'}=(P' \cap U) \setminus R_{P'}$.
    It is obvious that if $\dim P' \cap U =m-t+1$
    then $R_{\partial P'}$ is closed  in $M$,
    the images of the maps in
    ${\cal F}_{\partial P'}$  are contained in
    $\partial (P' \cap U) \setminus R_{ \partial P'}$ and form
     a discrete family in  $M \setminus R_{\partial P'}$.
     Now  assume that the last properties hold for a finite intersection $P'$ of $\cal P$
     with $\dim P' \cap U > m-t$.
    Then applying \ref{desc}  to all the maps of ${\cal F}_{ \partial P'}$
    we can enlarge $R_{\partial P'}$ to a closed  subset
      $R_{P'}$ of $M$ contained  in $P'\cap U$ and find extensions ${\cal F}_{P'}$ of  the maps
      in ${\cal F}_{\partial P'}$  such that
      $\dim R_{P'}\leq \dim(P' \cap U) -l+1=\dim(P' \cap U)+q-m-1$,
    $R_{P'}\setminus R_{\partial P'}$
    is a PL-subcomplex of $M \setminus R_{\partial P'}$,
      the images of the extensions ${\cal F}_{P'}$
    of the maps in ${\cal F}_{\partial P'}$  are contained in
    $G_{P'}$ and
    form a discrete family in  $M\setminus R_{ P'}$.
      Thus we can  define
    $R_P^{\Delta'}=R_{P'}$ for every finite intersection $P $ of $\cal P$ with
    $\dim P \cap U =m-t$ and   $\Delta' \subset \partial \Delta$
    such that $P(\Delta')=P'$.
    Then  proceeding  by induction on $\dim P'\cap U$  we  construct
    for every  intersection
    $P$ of $\cal P$ with $\dim P \cap U =m-t$ and $\Delta' \subset \partial \Delta$
    the collection of maps ${\cal F}_P^{\Delta'}$ and
   the set $R_P^{\Delta'}$
   such that
   $R^{\partial \Delta}$
   is a closed PL-presented subset of $M$ contained in $U$ such that
   $\dim R^{\partial \Delta}  \leq q-1$ and the images of the maps in
   ${\cal F}^{\partial \Delta}$
   are contained in $(M\cap U)\setminus R^{\partial \Delta}$ and
   form a discrete family in $M\setminus R^{\partial \Delta}$.
   From the construction it is clear that $R^{\partial \Delta}$
   is nowhere dense on the finite intersections of $\cal P$.

   By  \ref{desc} the set $R^{\partial \Delta}$ can be enlarged
   to a closed PL-presented subset ${R}$ of $M$ and  each map
   $f_{\partial \Delta}: S^{q-1} \lo M\setminus R^{\partial \Delta}$
   in ${\cal F}^{\partial \Delta}$ can be extended to
   a map $f'_{\Delta} : B^q \lo M\setminus {R}$
   such that $\dim R \leq q$, $R$ is
   nowhere dense on the finite intersections
   of ${\cal P}$ and the images of the maps $f'_{ \Delta}$ form a discrete
   family in $M\setminus {R}$.
    Hence for every $f_{\partial \Delta}\in {\cal F}^{\partial \Delta}$
     there is an open subset $Q(f_{\partial \Delta})$ of $M\setminus R$
     containing the image of $f_{\partial \Delta}$ such that
     ${\cal Q}=\{ Q(f_{\partial \Delta}):
     f_{\partial \Delta}\in {\cal F}^{\partial \Delta}\}$ is a discrete
     family in $M\setminus R$ and $f_{\partial \Delta}$
      is
   null-homotopic in $Q(f_{\partial \Delta})$.
   Then by \ref{hole1}  the black hole modification
     used in \ref{connect} and involving $f_{\partial \Delta}$
   can be carried out inside $Q(f_{\partial \Delta})$.

     Note that since $R$ is a PL-presented set of $\dim  \leq q$
   and $m \geq 2q+1$
   then for every finite intersection $P$ of $\cal P$ that
   intersects $U$ we have that
   the inclusion $(P\cap U) \setminus R \subset P\cap U$
   induces an isomorphism of the homotopy groups in
   $\dim \leq \dim (P \cap U) -(l-1)$ (=co-dimensions $\geq l-1$).

 Thus after removing $R$ from $M$, $F$ and  the elements of ${\cal P}$
    the black hole modification \ref{black}
     used in \ref{connect}
    can be applied independently
   for every map in ${\cal F}^{\partial \Delta}$
    in order to modify
   $\cal P$  to
   a decomposition ${\cal P}'$ of $M'=M \setminus R$
   and $F$ to a set $F'$
   such that ${\cal P}'=$the result of  all the modifications ${\cal P}^\pi$
   and
   $F'$=the result  of all the modifications $F^\pi$.
    Then it is easy to derive
   from \ref{connect} that
   $F'$ is a PL-subcomplex of $M'$
   of $\dim \leq\max\{m-q-1,\dim F\}$
   lying in $\Int M'$,
   ${\cal P}'$ is a partition on $M'\setminus F'$ and
   for every  intersection $P$ of
   $\cal P$ with $ P\cap U \neq \emptyset$,
    $P$  is modified to $P'$
  (= the result of all the modifications
   $P^\pi$)  such that
   $P'\cap (M'\setminus F')$ is  $l$-co-connected   and
  $P'\cap(M'\setminus F')$ is $(l-1)$-co-connected if
   $\dim P\cap U \geq m-t$.

   Note that
   the natural one-to-one correspondence between
   ${\cal P}$ and ${\cal   P}'$  defined by sending $P\in {\cal P}$
    to its  modification $P' \in{\cal P}'$  turns into a matching
   of partitions when ${\cal P}$ and ${\cal P}'$ are restricted to
   $M \setminus F$ and $M' \setminus F'$.

   It is clear that the images of the maps in ${\cal F}^{\partial \Delta}$
   are contained in the elements of $\st ({\cal P}, {\cal P})$.
   Assume that $\st ({\cal P}, {\cal P})$ refines an open cover ${\cal W}$
   of $M$ such that for every $A \in \st ({\cal P}, {\cal P})$ there is
   $W \in {\cal W}$ such that $A \subset W$ and the inclusion of $A$
   into $W$ induces the zero-homomorphism
   of the homotopy groups in $\dim \leq q-1$. By \ref{desc} the collection
   ${\cal Q}$ can be chosen so that ${\cal Q}$ refines ${\cal W}$.
   Then it is easy to see
   that
   ${\cal P}'$  refines $\st({\cal W}, {\cal W})$.

  \subsection{ Improving the total  connectivity
   of a partition }\label{total}
  Let $M$ be a triangulated $(q-1)$-connected $m$-dimensional manifold
  with $m \geq 2q+1$
  and let $l=m-q+2$. Assume that
  $F$ is a  PL-subcomplex of $M$ lying in
   $\Int M $
  and $\cal P$ is a decomposition of $M$
  such that  $\dim F  \leq m-q$ and
  $\cal P$ is an $l$-co-connected partition on
  $U=M \setminus F$.

  Apply  \ref{simul}  to improve the connectivity of the elements
  of ${\cal P}|U$ to the $(l-1)$-connectivity,  again apply
  \ref{simul} to the modified decomposition to improve the
  connectivity of the  intersections  of $\dim= m-1$
  (of the modified decomposition restricted to the complement
   of the modified $F$)
  to
  the $(l-1)$-connectivity and thus proceed on the dimension
  of the intersections $\geq l-1$  until we modify $M$ to
   an open subset $M'\subset M$,
   ${\cal P}$ to a decomposition ${\cal P}'$ of $M'$
  and  $F$ to
  a PL-subcomplex   $F'$ of $M'$ lying in
  $\Int M' $
  such that
  $M\setminus M'$ is a closed  PL-presented subset  of $M$,
   $\dim (M \setminus M')\leq q$,
    $\dim F' \leq \max  \{m-q-1, \dim F \}\leq m-q$,
   ${\cal P}'$   is
  an $(l-1)$-co-connected  partition on $U'=M' \setminus F'$
  and ${\cal P}$ admits a natural one-to-one correspondence
  to
  ${\cal P}'$ which sends
  each element of ${\cal P}$ to its modification in
  ${\cal P}'$ and which becomes a matching of partitions when
   ${\cal P}$ and ${\cal P}'$ are restricted to $U$
  and $U'$ respectively.

  \subsection{ Absorbing  simplexes}\label{absorb}

 Let  $M$ be a triangulated $(q-1)$-connected $m$-dimensional manifold
  with $m \geq 2q+1$
  and let $l=m-q+1$. Assume that
  $F$ is a  PL-subcomplex of $M$ lying
  in $\Int M $ such that $U=M \setminus F$
  is $l$-co-connected and  $\dim F  \leq m-q$ and assume that
   $\cal P$ is a decomposition of $M$
  such that
  $\cal P$ is an $l$-co-connected partition on $U$.
 Fix a triangulation $\cal T$ of $M$ for which
  $F$ is a subcomplex
         and let
          $F'$ be the $(m-q-1)$-skeleton of $F$.
          We will show how to modify $M$ to an open subset
  $M' \subset M$ and  ${\cal P}$ to a decomposition ${\cal P}'$
  of $M'$
   such that $M\setminus M'$ is a PL-subcomplex of $M$ of
  $\dim \leq q$,
   ${\cal P}'$ restricted to
  $U'=M' \setminus F'$
  is an
  $l$-co-connected partition and
  ${\cal P}$ admits a natural one-to-one
  correspondence to ${\cal P}'$ which becomes a matching of
   partitions when
  $\cal P$ and ${\cal P}'$ restricted to $U$ and $U'$ respectively.

  Assume that the triangulation $\cal T$ that we fixed in $M$ is the second
  barycentric subdivision of a triangulation for which $F$ and
  the elements of $\cal P$ are subcomplexes.
  Let
   ${\cal T}_F$ be the collection of all
   $(m-q)$-dimensional simplexes lying in $F$,
   $S_\Delta$ the link of $\Delta \in {\cal T}_F$ with respect to
   ${\cal T}$,
   $Q=$ the union of
  $\Int(\Delta *  S_\Delta)$ for all $\Delta $ in ${\cal T}_F$
  and $N = M\setminus (Q \cup F)$.
  Then
  for every $\Delta \in {\cal T}_F$  the link
  $S_\Delta$ of $\Delta$ is  an $(q-1)$-dimensional sphere  lying
  in $\Int M $,
  for different  simplexes
  $\Delta_1$ and
  $\Delta_2$
   in ${\cal T}_F$  the joins $\Delta_1 * S_{\Delta_1}$ and
  $\Delta_2 * S_{\Delta_2}$ do not  intersect on $U$,
   for every finite intersection $P$ of $\cal P$
  that intersects $U$ and for  every $\Delta \in {\cal T}_F$,
  $(P\cap U) \setminus \Int(\Delta * S_{\Delta})$
   is a deformation
  retract of $P \cap U$. Note that
    $N$ is a manifold which is a deformation retract of $U$ and therefore $N$
  is $l$-co-connected.

   Fix $\Delta \in {\cal T}_F$ and let an open subset $V_\Delta$
   of $N$ be such that $S^{q-1}=S_\Delta \subset V_\Delta$
   and the inclusion of $S^{q-1}$ into
   $V_\Delta$ is null-homotopic in $V_\Delta$. Then by \ref{hole2}
   there are an element $P_\Delta \in \cal P$ and a PL-embedding
   of a cube $B^m=B^q \times B^{m-q}$ into $\Int M$
     such that $B^m \subset (\Delta * S_\Delta) \cup V_\Delta$,
   $\Delta =0 \times B^{m-q}=F \cap B^m$, $P \cap B^m=(P \cap B^q)\times B^{m-q}$
   for every $P \in \cal P$,
   $\partial B^q \times B^{m-q} \subset \Int P_\Delta$,
   $\cal P$ restricted to $(M \setminus \Int B^m)\cap U$
   and $\partial B^m \cap U$ are partitions and,
   $(P \setminus \Int B^m) \cap U$ is
   a deformation retract of $(P \setminus \Int B^q ) \cap U$
   for every finite intersection $P$ of  $\cal P$.

   Let us briefly describe the general idea of the construction in \ref{hole2}.
    The inclusion
   of $S^{q-1}$ can be extended to a PL-embedding of
   a $q$-dimensional   cube $ B_\#^q$ into $V_\Delta$ such that
   $S^{q-1}=\partial B_\#^q$ and the embedding of $B_\#^q$
   extends to a PL-embedding of
   an $m$-dimensional cube  $B_\#^m=B_\#^q \times B_\#^{m-q}$
   into $V_\Delta$  such that $\partial B_\#^q \times B_\#^{m-q}
   \subset B_\#^m \cap \partial (\Delta * S_\Delta)$ and
   $P\cap B_\#^m =(P \cap B_\#^q) \times B_\#^{m-q}$ for every
   $P \in \cal P$. Let $P_\Delta \in \cal P$ be such that
   $\Int (P_\Delta  \cap B_\#^q) \neq \emptyset$
   and let a $q$-dimensional cube $\sigma_{\#}$ be PL-embedded in
   $\Int (P_\Delta  \cap B_\#^q)$. Then the cube
   $B^m=(\Delta * S_\Delta)\cup ((B_\#^q \setminus
   \Int \sigma_{\#}) \times B_\#^{m-q})$
    can be represented as the product
    $B^m=B^q \times B^{m-q}$ with the required properties.
    See \ref{hole2} for details.

    Let $B^{m-q}_*$  be  a cube lying in the geometric interior of
    $B^{m-q}$. Use the notations of \ref{black}.
    Define  $F_L=F \setminus \Int B^m$, $P_L = P\setminus \Int B^m$
    if $P \neq P_\Delta$ and $P \in {\cal P}$,
    $(P_\Delta)_L=(P_\Delta \setminus \Int B^m)\cup
     (B^q \times
    B^{m-q}_*)$
    and define
    ${\cal P}_L=\{P_L : P \in {\cal P} \}$. Then
    $S=\partial B^{m-q}=\partial \Delta \subset F_L$
     and it is easy to see that
     ${\cal P}_L$ is
    an $l$-co-connected partition on $L\setminus F_L$ and
    ${\cal P}_L$ is a partition on $\partial T \setminus F_L$.
    Now apply the black hole modification with
    the black hole $S=\partial \Delta$  and denote $F^\pi =F^\pi_L$,
    $P^\pi =P_L^\pi$ for $P \in \cal P$ and
    ${\cal P}^\pi =\{P^\pi : P \in {\cal P}\}= {\cal P}_L^\pi$. Note
    that $F^\pi_L =F_L$ since $F \cap \partial B^m=\partial \Delta$.
    Then  ${\cal P}^\pi$ is a decomposition of $M$ which
    is an $l$-co-connected partition on the manifold
    $M\setminus F^\pi=(M \setminus F)\cup \Int \Delta$
   in which  we have
    absorbed the geometric interior of $\Delta$ in the
    modification of $P_\Delta$. Note that the one-to-one
    correspondence between ${\cal P}$ and ${\cal P}_L$ defined
    by $P\lo P_L$, $P \in \cal P$ becomes a matching of
    partitions when ${\cal P}$ and ${\cal P}_L$ are restricted
    to $M \setminus F$ and $L \setminus F_L$ respectively.
    Then by \ref{black} the one-to-one
    correspondence between ${\cal P}$ and ${\cal P}^\pi$ defined
    by $P\lo P^\pi$, $P \in \cal P$ becomes a matching of
    partitions when ${\cal P}$ and ${\cal P}^\pi$ are restricted
    to $M \setminus F$ and $M \setminus F^\pi$ respectively.

    Note that each $S_\Delta$ can be contracted to a point in $N$
    outside some  neighborhood of $F$ in $M$.
    Then by \ref{desc} there is a PL-subcomplex  $R$ of $M$ contained in $\Int N$
    and
    nowhere dense on the finite intersections of $\cal P$ with
    $\dim R \leq q$
     such that
    the collection ${\cal V}$ of the sets $V_\Delta, \Delta \in {\cal T}_F$
    can be chosen to be
    discrete in $N\setminus R$. Then
  replacing  $M$ by $M'=M \setminus R$ and removing $R$ the elements
  of $\cal P$ we can
   perform all the black hole modifications
    independently for all $\Delta \in {\cal T}_F$
     and    get from all the modifications   ${\cal P}^\pi$
    the corresponding  decomposition ${\cal P}'$ of $M'$ with
    the required properties.
     It is clear that the natural one-to-one correspondence
     between $\cal P$ and ${\cal P}'$ defined by sending each element of
     $\cal P$ to its modification  becomes
    a matching of partitions  when the decompositions are restricted to
    $M\setminus F$ and $M'\setminus F'$.

    Let   $\cal W$  be an open cover of   $M$ such that
    $\st({\cal P, {\cal P})}$ refines ${\cal W}$  and assume that
    for every $S_\Delta$, $\Delta \in {\cal T}_F$ there is
    $W_\Delta \in \cal W$ such that $S_\Delta \subset W_\Delta$
    and $S_\Delta$ can be contracted to a point in
    $W_\Delta \cap N$. Then the collection ${\cal V}$
    can be chosen so that $\cal V$ refines $\cal W$
    and this implies that ${\cal P}'$  refines
    $\st ({\cal W}, {\cal W})$.

   \subsection{ Improving connectivity via  a matching }\label{twin}

 Let $M_i$  be an $l_i$-co-connected triangulated manifolds
 such that
  $m_i =\dim M_i \geq 2(m_i -l_i)+3$, $i=1,2$ and $m_1-l_1=m_2-l_2$.
   Suppose that
   ${\cal P}_1$
  and ${\cal P}_2$ are  partitions of $M_1$ and $M_2$ respectively
   such that there is a matching between
   ${\cal P}_1$
  and ${\cal P}_2$, and
  ${\cal P}_1$ is $l_1$-co-connected. We  will show how to modify
  $M_2$ to $M'_2$ and ${\cal P}_2$ to ${\cal P}'_2$ such that
  $M'_2$ is an open $l_2$-co-connected submanifold of $M_2$,
  $M_2\setminus M'_2$ is a PL-presented
  closed subset of $M_2$,
  $\dim (M_2 \setminus M'_2) \leq l_2-2$
   and
  ${\cal  P}'_2$ is an $l_2$-co-connected partition
  of $M'_2$ which admits a natural  matching  to
  $ {\cal P}_2$ defined
  by sending each element of ${\cal P}_2$  to
  its modification in ${\cal P}'_2$.

  By \ref{total} modify $M_2$ to $M_2^0 \subset M_2$, ${\cal P}_2$
  to ${\cal P}_2^0$  and  construct a subset $F^0$ of $M_2^0$
  such that $M_2\setminus M_2^0$ is
   a  PL-subcomplex of $M_2$,
   $\dim (M_2\setminus M_2^0)) \leq 1$,
  $F^0$ is a PL-subcomplex  of $M_2^0$,
  $\dim F^0 \leq m_2 -2$,
    ${\cal P}_2^0$ is
   a decomposition of $M_2^0$ such that
   ${\cal P}_2^0$ restricted to
     $M_2^0 \setminus F^0$ is an $m_2$-co-connected
     partition
   that admits a matching to ${\cal P}_2$ and therefore
   to ${\cal P}_1$.

   Now assume that for $ 0 \leq t < m_2-l_2$  we have
   constructed $M_2^{2t}$, ${\cal P}_2^{2t}$ and $F^{2t}$
   such that
   that $M_2\setminus M_2^{2t}$ is
   a  PL-presented closed subset of $M_2$,
   $\dim (M_2\setminus M_2^{2t}) \leq t+1$,
  $F^{2t}$ is a PL-subcomplex  of $M_2^{2t}$,
  $\dim F^{2t} \leq m_2 -t-2$,
    ${\cal P}_2^{2t}$ is
   a decomposition of $M_2^{2t}$ such that
   ${\cal P}_2^{2t}$ restricted to
     $M_2^{2t} \setminus F^{2t}$ is an $(m_2-t)$-co-connected
     partition
   that admits a matching to ${\cal P}_2$ and therefore
   to ${\cal P}_1$.  Proceed to $t+1$ as follows.

  By \ref{total} modify $M_2^{2t}$ to
  $M_2^{2t+1} \subset M_2^{2t}$, ${\cal P}_2^{2t}$
  to ${\cal P}_2^{2t+1}$, $F^{2t}$ to $F^{2t+1}$
  such that $M_2^{2t}\setminus M_2^{2t+1}$ is
   a  PL-subcomplex of $M_2^{2t}$,
   $\dim (M_2^{2t}\setminus M_2^{2t+1}) \leq t+2$,
  $F^{2t+1}$ is a PL-subcomplex  of $M_2^{2t+1}$,
  $\dim F^{2t+1} \leq m_2 -t-2$,
    ${\cal P}_2^{2t+1}$ is
   a decomposition of $M_2^{2t+1}$,
   ${\cal P}_2^{2t+1}$ restricted to
     $M_2^{2t+1} \setminus F^{2t+1}$ is an $(m_2-t-1)$-co-connected
     partition
   that admits a matching to ${\cal P}_2^{2t}$ restricted
   to $M_2^{2t} \setminus F^{2t}$ and therefore
   to ${\cal P}_2$ and ${\cal P}_1$.

  Then since   ${\cal P}_1$ is $l_1$-co-connected on the $l_1$-co-connected
     manifold $M_1$
     we conclude by \ref{match} that $M_2^{2t+1} \setminus F^{2t+1}$
     is $(m_2 -t-1)$-co-connected. Therefore, by \ref{absorb},
     $M_2^{2t+1}$, $F^{2t+1}$ and ${\cal P}_2^{2t+1}$ can be modified
     to $M_2^{2t+2}$, $F^{2t+2}$ and ${\cal P}_2^{2t+2}$ respectively
     such that  $M_2^{2t+2} \subset M_2^{2t+1}$,
    $M_2^{2t+1}\setminus M_2^{2t+2}$ is
   a  PL-subcomplex of $M_2^{2t+1}$,
   $\dim (M_2^{2t+1}\setminus M_2^{2t+2}) \leq t+2$,
  $F^{2t+2}$ is a PL-subcomplex  of $M_2^{2t+2}$,
  $\dim F^{2t} \leq m_2-t -3$,
    ${\cal P}_2^{2t+2}$ is
   a decomposition of $M_2^{2t+2}$,
   ${\cal P}_2^{2t+2}$ restricted to
     $M_2^{2t+2} \setminus F^{2t+2}$ is an $(m_2-t-1)$-co-connected
     partition
   that admits a matching to ${\cal P}_2^{2t+1}$
   restricted  $M_2^{2t+1} \setminus F^{2t+1}$
   and therefore
   to ${\cal P}_2$  and ${\cal P}_1$.

  Then   for $t=m_2-l_2$, $M'_2=M_2^{2t}\setminus F^{2t}$ and
   ${\cal P}'_2= {\cal P}_2^{2t}$ restricted to $M'_2$
   will have the required properties (note that
   $  m_2-l_2+1  \leq l_2-2$ and therefore
   $\dim (M_2\setminus M'_2) \leq l_2-2$).

     \subsection{ Moving to a rational position }\label{rational}

     Let $M$ be  a triangulated  manifold with the rational structure
     determined by a triangulation  $\cal T$ (see  Introduction).
     The PL-notions defined in \ref{general} can be translated to
     the corresponding rational notions  by referring to the rationsl
     structure instead of  the PL-structure.
      For example, $P\subset M$ is  a rational
      subspace if $P \subset M$ is rationally embedded in $M$,
       $P$ is a rational subcomplex
     if there is a rational triangulation of $M$ for which
     $P$ is a subcomplex, $\cal P$ is a rational decomposition of
     $M$ if the elements of $\cal P$ are rational subcomplexes and so on.

     Assume that  ${\cal T}'$ is  a  (not necessarily rational) triangulation
     of $M$.
     We will show that the identity map of $M$ can be
     arbitrarily closely approximated by  a PL-homeomorphism $f : M \lo M$
     such that
     for very $\Delta' \in {\cal T}'$, $f (\Delta')$ is a rational
     subcomplex of $M$.

          Embed $M$ into a Hilbert space by a  map
       which is linear on every simplex of  ${\cal T}$ and
        refer to this Hilbert space when
       properties of linearity are used.
     Let
     ${\cal T}''$ be a triangulation of $M$ such that ${\cal T}''$
     is a subdivision of both ${\cal T}$ and ${\cal T}'$, and
     the simplexes of ${\cal T}''$ are linear.
    Every vertex $v$ of ${\cal T}''$  approximate by
     a rational point $p_v$ (= a point with rational
      barycentric coordinates with respect to $\cal T$) such that for every
      simplex $\Delta$ of ${\cal T}$, $p_v \in \Delta$ if and only if
      $v \in \Delta$.
      The approximation of the vertices of ${\cal T}''$
      can be  extended to the PL-map $f : M \lo M$ sending each vertex
      $v$ of ${\cal T}''$ to $p_v$ such that $f$ is linear
      on each simplex of ${\cal T}''$. If $p_v$ is sufficiently close to $v$
      for every vertex $v \in {\cal T}''$ then the map $f$ is
      a PL-homeomorphism that can be chosen to be
      arbitrarily close to the identity map. Clearly $f$
       sends every simplex of ${\cal T}''$ to
      a rational simplex  and therefore $f$ sends every simplex of ${\cal T}'$
      to a rational subcomplex of $M$.

      Let   $M$ be an open subset of a space $Y$ and    $d$ a metric
       on  $Y$. We can assume that the homeomorphism $f: M \lo M$
       constructed above satisfies
       $d(y, f(y)) \leq d(y, Y\setminus M)$ for every $y \in M$. Then
       $f$ extends to the homeomorphism $g : Y \lo Y$ such that
       $g(y)=y$ if $y \in Y \setminus M$ and $g(y)=f(y)$
       if $y \in M$.
       Now assume that $R$ is a PL-presented subset of $M$ and let
      $R_1 \subset R_2 \subset \dots \subset R_n$
      be  closed subsets
      of $M$ such that
       $R=R_n$, $R_1$ is a PL-subcomplex of $M$ and $R_{i+1} \setminus R_i$
       is a PL-subcomplex of $M \setminus R_i$, $i=1, \dots, n-1$.
       Let us show  that
       the identity map of $Y$ can be arbitrarily closely
       approximated by
     a homeomorphism $g : Y \lo Y$ such that $g(y)=y $ for every $y \in Y\setminus M$,
     $g(R)$ is
     a closed  rationally presented subset of $M$ and $g$ restricted to $M\setminus R$
     is a PL-homeomorphism to $M\setminus g(R)$
     sending $\cal P$ to a rational decomposition
     of $M\setminus g(R)$ (note that the rational structure of $M$  induces
     the corresponding     rational structure on open subsets of $M$
     such that the inclusions are rational maps).

     Approximate the identity map of $Y$ by a homeomorphism $g_1 : Y \lo Y$
     such that $g_1$ does not move the points of $Y\setminus M$,
     $g_1$ restricted to $M$ is a PL-homeomorphism and
     $g_1(R_1)$ is a rational subcomplex of $M$. Approximate $g_1$ by
     a map $g_2 : Y\lo Y$ such that $g_2$ does not move  the points
     of $(Y \setminus M) \cup g_1(R_1)$, $g_2$ restricted to $M\setminus g_1(R_1)$
     is a PL-homeomorphism  and $g_2(g_1(R_2\setminus R_1))$ is a rational subcomplex
     of $M \setminus g_1(R_1)$. Proceed by induction and construct
     for every $i=1,2,\dots , n-1$ an approximation of  the  identity map
     of $Y$ by a   homeomorphism
     $g_{i+1} : Y \lo Y$ such that such that $g_{i+1}$ does not move  the points
     of $(Y \setminus M) \cup (g_i\circ \dots \circ g_1)(R_i)$,
     $g_{i+1}$ restricted to $M\setminus (g_i\circ \dots \circ g_1)(R_i)$
     is a PL-homeomorphism  and $(g_{i+1}\circ \dots \circ g_1)(R_{i+1}\setminus R_i)$
     is a rational subcomplex
     of $M \setminus (g_i\circ \dots \circ g_1)(R_i)$.
     And finally construct $g_n$ with the additional property that
      $g=g_n \circ \dots \circ g_1$
     sends ${\cal P}$ to a rational decomposition of $M \setminus g(R)$.
     Then $g$ will have the required properties.

  \end{section}
  \begin{section}{Proof of Theorem \ref{t2}}
  Let $A$ be a closed subset of a space $X$.
  A collection ${\cal C}$ of subsets $X\setminus A$ is said to properly
  approach $A$  if
  for every
  sequence $C_j$ of elements of ${\cal C}$ such that
   there is a sequence of points $x_j \in C_j $ converging
 to a point $a \in A$
  we have that $\lim_{j \rightarrow \infty} C_j =a$
  (that is $\lim_{j \rightarrow \infty} y_j =a$ for every
  sequence $\{ y_j\}$ with $y_j \in C_j$).
  Note that if $\cal C$ properly approaches $A$  then
  the closure in $X$  of each $C \in \cal C$ does not intersect
   $A$.
  Also  note that if collections  ${\cal C}$ and ${\cal C}'$ of subsets
  of $X \setminus A$
   properly approach $A$ then $\st({\cal C}, {\cal C}')$ properly
  approaches $A$ as well.

   \begin{proposition}
  \label{p2}
  Let $X$ be an $n$-dimensional \noo space,
   let $A$ be a $Z$-set in $X$  and let $\cal C$
   be a cover of $X \setminus A$ that properly approaches $A$.
   Then for every $C \in \cal C$  there is an open set
   $C\subset V_C \subset X \setminus A$ such that
   the inclusion $C \subset V_C$ induces the zero-homomorphism
   of the homotopy groups in $\dim \leq n-1$ and
    ${\cal V} =\{ V_C : C \in {\cal C} \}$ properly approaches $A$.
   \end{proposition}
   {\bf Proof.} Fix a metric $d$ on $X$ such that $d(x,y) \leq 1$
   for every $ x, y \in X$.
   If  $C \in \cal  C$ is not a singleton
   denote $d_C = \inf $ of ${\rm diam} G $ for
    open subsets $G $ of  $X$ such that  $C \subset G$ and the inclusion
   $C \subset G$ induces the zero-homomorphism of the homotopy
   groups in $\dim \leq n-1$ and if $C \in \cal  C$ is a singleton
   define $d_C =d(C,A)$.  For every
   $C \in \cal C$ fix an open set $G_C$ such that ${\rm diam} G_C \leq  2d_C$
   and the inclusion $C \subset G_C$ induces the zero-homomorphism
   of the homotopy groups in $\dim \leq n-1$.

   Let ${\cal U}$ be an open cover of $X \setminus A$ which
   properly approaches  $A$.
   Denote ${\cal C}_n = \{ C \in {\cal C}: 1/n+1 < d_C \leq 1/n \}$.
   Then the closure of the union of the elements
   of ${\cal C}_n$ does not intersect $A$.
   For a given $n$ approximate the identity map of $X$ by a closed
   embedding  $e_n : X \lo X$ such that $e_n(X) \subset X\setminus
   A$ and
   for  $C \in {\cal C}_n$,
   ${\rm diam} (e_n(G_C)) \leq {\rm diam} G_C
   +1/n$, $e_n(C) \subset \st(C, {\cal U})$
    and  every map $f : S^p \lo C$
    from a $p$-dimensional sphere $S^p$, $p \leq n-1$  can be homotoped into
   $e_n(C)$ inside $\st(C, {\cal U})$.
      Denote
   $Y_C = \st(C, {\cal U}) \cup e_n(G_C)$ for $C \in {\cal C}_n$.

   Then the inclusion $C \subset Y_C$
   induces
   the zero-homomorphism of the homotopy groups in $\dim \leq n-1$
   and the cover $\{ Y_C : C \in {\cal C} \}$ of $X \setminus A$
    properly
   approaches $A$. Enlarge  the sets $Y_C$   to open
   sets $V_C$ such that ${\cal V}= \{ V_C : C \in {\cal C} \}$
   properly approaches $A$. Then ${\cal V}$ has the required
   properties.
      \hfill $\Box$
  \\\\

  Suppose that  $A$ is a closed subset of $X$ and  $X \setminus A$ is
  embedded as a dense subset in $M$. Let us show that there is an open subset $V$
  of $M$ containing $X \setminus A$  for which there is  a space $Y$ such that
    $X$ and $V$ embed into
    $Y$ such that  $Y=X \cup V$,  $X\setminus A = X \cap V$, $A=Y \setminus V$ and
    $A$ is
  closed in $Y$.

  Embed $X$ into a Hilbert space and let   a metric $d(x,y)=|x-y|$
  on $X$ comes from a norm in $H$.
  For every $x \in X \setminus A$ define
  $W(x)=\{y \in X : d(x,y) < d(x,A)/2\}$
  and let ${\cal W} =\{ W(x) : x \in X \setminus A\}$.
  Refine ${\cal W}$ by a locally finite open cover
  $\cal G$ of $X \setminus A$.
  For every $G \in \cal G$ pick up a point $x_G \in X\setminus A$
     such that
    $G \subset W(x_G)$ and observe that
    if $x$ belongs to  $ G \in \cal G$
    then $ d(x_G,A)/2 \leq d(x,A) \leq 3d(x_G,A)/2$ and
    $d(x_G, x)\leq d(x,A)$.
    Extend each $G \in \cal G$ to an open subset $V_G $ of
    $M$ such that $G = (X \setminus A) \cap V_G$
    and ${\cal V}= \{ V_G : G \in {\cal G}\}$
    is  locally finite  in  $V =\cup \{ V_G : G \in {\cal G}\}$.
   Let $\{ f_G : G \in {\cal G}\}$ be a partition of unity
   subordinated to $\cal V$.
     Set
   $p_G=(x_G, d(x_G, A)) \in H \times (0,\infty)$ and
   define a map $\phi_V : V \lo H\times (0, \infty)$
   by $\phi_V (x)=\sum_{G \in {\cal G}} f_G(x)p_G$.
   Identify $H$ with $H\times 0 \subset H \times [0,\infty)$
   and let
    $\cal U$ be an open cover of $H \times (0, \infty)$
   that properly approaches $H$. Approximate $\phi_V$ by
   a $\cal U$-close embedding $h_V : V \lo H\times(0, \infty)$.
   Define a set  $Y$ as the free union of $A$ and $V$. Then $X$ is a subset
   of $Y$.
   Denote by  $h_A$  the inclusion of $A$ into $H$
   and define a function  $h : Y \lo H\times [0, \infty)$ such that
   $h|A=h_A$ and $h|V=h_V$.

   Let us show that  $h|X$ is an embedding of $X$ into
   $H \times [0, \infty)$.
   To show this we need to verify that for a sequence
    $ x_j $  in
   $ X\setminus A$
   and $a \in A$,
   $\lim x_j =a$ if and only if $\lim h_V(x_j)=a$.
   Since
   $\cal U$ properly approaches $H$, $\lim h_V(x_j)=a$ is
   equivalent to $\lim \phi_V (x_j)=a$.
   For every $G \in \cal G$
   containing $x_j$ we have that $d(x_j, A)/2 \leq d(x_G, A)\leq 3d(x_j,A)/2$
   and  $d(x_G,x_j) \leq d(x_j,A)$. Then
      $ d(x_j, A)/2  \leq  \sum_G f_G(x_j)d(x_G, A) \leq 3d(x_j,A)/2$
      and $|\sum_G f_G(x_j)x_G - x_j| \leq d(x_j, A)$.
    Assume that $\lim x_j=a$.
    Then $\lim d(x_j, A)=0$ and hence
 $\lim_j \sum_G f_G(x_j)d(x_G, A)=0$ and
     $\lim_j \sum_G f_G(x_j)x_G=\lim x_j=a$.
    Thus $\lim \phi_V (x_j)=a$.
     Now assume $\lim \phi_V (x_j)=a$. Then
  $\lim_j \sum_G f_G(x_j)d(x_G, A)=0$
  and $\lim_j \sum_G f_G(x_j)x_G=a$. This implies  $\lim d(x_j, A)=0$
  and   $\lim_j \sum_G f_G(x_j)x_G=\lim x_j=a$. Thus we showed that
  $h$ embeds $X$ into $H \times [0,\infty)$.
  Now we can transfer the topology of $h(Y)$ to $Y$ and  get
    the required space.\\

  Let $M_i$  be an $l_i$-co-connected triangulated manifold, $i=1,2$
 such that  $m_i \geq 2(m_i -l_i) +3$
 and $m_1-l_1 =m_2-l_2$
 for $m_i =\dim M_i$.
 Denote $k_i=l_i-2$ and  $n = m_i-l_i +1$.
 Then by Theorem \ref{t1}, $M_i(k_i)$ is an $n$-dimensional
 \no space. Assume that  $X_i$ and $A_i$ satisfy the assumptions of Theorem \ref{t2},
  let $f: A_1 \lo A_2$ be
 a homeomorphism   and let  $X_i \setminus  A_i$ be
 homeomorphic to $M_i(k_i)$.  Identify $X_i\setminus A_i$ with $M_i(k_i)$ and,
  by the property proved above,   replace $M_i$ by an open
 subset of $M_i$ containing $M_i(k_i)$ and  assume that $X_i$ and $M_i$
 are subspaces of a space $Y_i$ such that
 $Y_i = X_i \cup M_i$, $X_i\setminus A_i = X_i \cap M_i$,
  $A_i=Y_i \setminus M_i$ and $A_i$  is closed in $Y_i$.

 Let $\mu : {\cal C}_1 \lo {\cal C}_2$ be
  a one-to-one correspondence
  of
 between collections of subsets of $M_1$ and $M_2$ respectively.
 We say that $\mu$
 agrees with $f : A_1 \lo A_2$ if for every sequence $C_j, j=1,2,\dots$
 of elements of ${\cal C}_1$,
  $C_j$ converge to $a \in A_1$ in $Y_1$  (as $j \rightarrow\infty$) if and only
  $\mu(C_j)$ converge to $f(a)$ in $Y_2$. It is easy to check
   that if ${\cal C}_1$ and
  ${\cal C}_2$ properly approach $A_1$ and $A_2$ respectively,
  $\mu$ agrees with $f$ and a collection ${\cal C}$ of subsets of
  $M_1$ properly approaches $A_1$ then $\mu(\st({\cal C}, {\cal C}_1))$
  properly approaches $A_2$ where
  $\mu(\st({\cal C}, {\cal C}_1))=\{ \mu(\st(C,  {\cal C}_1)) : C \in {\cal C} \}$
  and $\mu(\st(C,  {\cal C}_1))=\cup \{\mu(C'):
  C'\in {\cal C}_1, C'\cap C \neq \emptyset \}$.

   Let $M$ be a triangulated space. It follows from the interpretation
   of $M(k)$ given in Introduction
   that if $P$  is either a rational subcomplex or an open subset of
  a triangulated space $M$
   then $P \cap M(k)=P(k)$
  where $P$ is considered with the rational structure for which
  the inclusion of $P$ into $M$ is a rational map
  (as  a matter of fact this is true for every rational subspace $P$
  of $M$ however we will not use this fact).

 The poof of Theorem \ref{t2} is based on
 \begin{proposition}
 \label{p1} Let  $M_i, A_i, X_i,Y_i, k_i, l_i, m_i $ and $f$  be as above and
 let ${\cal P}_1$
 be an $l_1$-co-connected rational partition of an
 open submanifold $M'_1$ of $M_1$ such that $M_1(k_1) \subset
 M'_1$, ${\cal P}_1$ has no non-empty intersections of
 $\dim \leq l_1 -2$
 and ${\cal P}_1$ properly approaches  $A_1$. Then there are an
 open submanifold $M'_2$ of $M_2$ and an $l_2$-co-connected
 rational partition of $M'_2$ such that $M_2(k_2) \subset M'_2$,
 ${\cal P}_2$ properly approaches $A_2$ and ${\cal P}_2$
 admits a matching  $\mu : {\cal P}_1 \lo {\cal P}_2$ which agrees
 with $f: A_1 \lo A_2$.
 \end{proposition}

 Let us show how Theorem \ref{t2} can be derived from Proposition
  \ref{p1}.
 \\
  {\bf Proof of Theorem \ref{t2}.}
  Assume that theorem holds in $\dim \leq n-1$ and let us prove it in
  $\dim=n$.
  Fix complete metrics on $X_1$ and $X_2$ and let $\epsilon >0$.
  Take a locally finite cover $\cal C$ of $X_1 \setminus A_1$ of
   mesh${\cal C} < \epsilon$
  such that $\cal C$
   properly approaches $A_1$. Then one can find an open subset
  $M'_1$  of  $M_1$ such that $M_1(k_1) \subset M'_1$ and the
  closures
   in $M'_1$ of the elements of $\cal C$ form a locally finite cover
   of  $M'_1$. Now there is a rational triangulation  of $M'_1$
   such that  the decomposition of $M'_1$ formed from the $m_1$-dimensional simplexes
   of the triangulation  properly approaches $A_1$ and
   the diameter of the intersection of every $m_1$-dimensional simplex
   with $M_1(k_1)$ is $< 2\epsilon$ with respect
   to the complete metric fixed in $X_1$. Define ${\cal P}_1$ as the stars
   of the vertices of this triangulation with respect to its first
   barycentric subdivision. Then ${\cal P}_1$ is a
   rational  partition of $M'_1$ which properly approaches $A_1$.
   Note that  the intersections of ${\cal P}_1$ of $\dim \leq l_1-2$
   do not intersect $M_1(k_1)$ and therefore
   by removing them  from $M'_1$ and from the elements of ${\cal P}_1$
   we may assume
   that ${\cal P}_1$ has no non-empty intersections
   of $\dim \leq l_1-2$.
   Note all the intersections of ${\cal P}_1$ remain contractible
   and hence ${\cal P}_1$ is $l_1$-co-connected.

   By Proposition \ref{p1}
  there are an
  open submanifold $M'_2$ of $M_2$ and an $l_2$-co-connected
  rational partition of $M'_2$ such that $M_2(k_2) \subset M'_2$,
  ${\cal P}_2$ properly approaches $A_2$ and ${\cal P}_2$
  admits a matching  $\mu : {\cal P}_1 \lo {\cal P}_2$ which agrees
  with $f: A_1 \lo A_2$.

  We are going to  construct  homeomorphisms
  $f_P :  P(k_1) \lo \mu(P)(k_2)$
   for  all finite intersections
  $P$  of ${\cal P}_1$ of $\dim \leq m_1 -1$
  which will agree on the common intersections.
  Let $t \leq  m_1-1$ and assume that for every
  finite intersection $P$ of $\dim \leq t-1$
   we already constructed $f_P$. Take an intersection
   $P$ of ${\cal P}_1$  such that $\dim P=t$.
   By Theorem \ref{t1},  $P(k_1)$ is a \no space
   of $\dim \leq n-1$. Note  for the union $P'$  of the
   intersections
   of ${\cal P}_1$ of $\dim < t$
   that are contained in $P$,   $P'$ lies in $\partial P$ and
   hence $P'(k_1)$ is a $Z$-set in $P(k_1)$.
   Define the homeomorphism  $f_{P'}: P'(k_1) \lo \mu(P')(k_2)$  by the
   homeomorphisms of the intersections forming $P'$.
      According to our assumption Theorem \ref{t2} holds
      in $\dim \leq n-1$. Therefore $f_{P'}$ can be
      extended to a homeomorphism  $f_P : P(k_1) \lo \mu(P)(k_2)$.

      Recall that ${\cal P}_1$ and ${\cal P}_2$ properly
      approach $A_1$ and $A_2$ respectively and the matching $\mu$
      agrees with $f$.  Then a homeomorphism extending $f$ between $X_1$
       and  $X_2$ can be obtain by pasting homeomorphisms
      from   $ P(k_1)$ to $ \mu(P)(k_2)$ for $P \in {\cal P}_1$
       which extend already defined  homeomorphisms on intersections of
       ${\cal P}_1$ of $\dim \leq m_1-1$.

       Fix $P \in {\cal P}_1$ and let $P'$ be the union of
       the intersections of ${\cal P}_1$ of $\dim \leq m_1-1$ that
       are contained in $P$.
       The homeomorphism  $f_{P'} : P'(k_1) \lo \mu(P')(k_2)$
       is a homeomorphism of $Z$-subsets of $P(k_1)$
       and $\mu(P)(k_2)$ respectively. Therefore we can
       repeat for  $P$ and $\mu(P)$ the same procedure  that we did for $M_1$ and
       $M_2$ but this time in the opposite direction from $\mu(P)$ to $P$
       first "splitting" $\mu(P)$ into small pieces   and then defining
       the corresponding splitting of $P$.

       Thus going back and force we are able after each iteration to extend $f$
       to  a partial homeomorphism  of  bigger and bigger
       parts of $X_1$ and $X_2$ and  simultaneously
       to restrict for each point of $X_1$ and $X_2$ the
       set to which this point can be sent under a possible extension of
       $f$.  Finally, passing to
       the limit and using the completeness of $X_1$ and $X_2$
       we get the desired homeomorphism.
        \hfill $\Box$
  \\\\
  {\bf Proof of Proposition \ref{p1}}. The proof of the
  proposition splits into two independent parts.
  \\\\
  {\bf  Constructing  an initial partition ${\cal P}_2$}.
  Let $P_0, \dots , P_t$ be distinct elements of ${\cal P}_1$ and
  let $P=P_0 \cap \dots \cap P_t$. Denote by  $\Delta_P$  a
  $t$-dimensional simplex with vertices $ v_0, \dots, v_t$ and
  for a simplex $\Delta'\subset \Delta_P$ denote by $V_{\Delta'}$ the set
  of vertices of $\Delta'$.
  A map
  $e_P : \Delta_P \lo  M'_1$
  is  called  a map witnessing the intersection $P$ if
  $e_P(\Delta_P) \subset \Int(P_0 \cup \dots \cup P_t)$,
  $e_P(v_i) \in \Int P_i$ and
  $P_i$ does not intersect $e(\Delta^i_P)$ where
  $\Delta^i_P$ is the face of $\Delta_P$  spanned by
  the vertices $\{ v_0, \dots, v_t \} \setminus \{ v_i \}$,
  $i=0,\dots , t$. Note that if there is a map  $e_P$ witnessing
   the intersection $P$ then  $P \neq \emptyset$.
   Indeed,
   aiming at a contradiction assume that $P=\emptyset$.
   Enlarge each $e^{-1}(P_i)$ to an open subset $e^{-1}(P_i)\subset G_i$
   of $\Delta_P$ such that $G_i$ does not meet $\Delta^i_P$ and
   $G_0 \cap \dots \cap G_t =\emptyset$.
   Let $f_0,f_1, \dots, f_t$ be a partition of unity
   subordinated to the cover ${\cal G}=\{ G_0, \dots G_t\}$ of $\Delta_P$.
   Define $f : \Delta_P \lo \Delta_P$ by
   $f(x)=f_0(x)v_0+\dots+f_t(x)v_t$, $x \in \Delta_P$. Then
   $f(\Delta_P)\subset \partial \Delta_P$ and, for every $i$
   and $x \in \Delta^i_P$, $f(x) \in \Delta^i_P$. Hence
   $f$ restricted to $\partial \Delta_P$ is homotopic to
   the identity map of $\partial \Delta_P$ and this contradicts
  the fact    that $\partial \Delta_P$ is not a retract of $\Delta_P$.

  Assume that $P$ is a non-empty intersection and let us show how
  to construct a map witnessing the intersection $P$.
  Let $F_i \subset  \Delta_P$   be the star of $v_i$ with
  respect to the first barycentric subdivision of $\Delta_P$.
  For a simplex $\Delta' \subset \Delta_P$ denote
  $C_{\Delta'}=\cap \{F_i : v_i \in V_{\Delta'}\}\subset \Delta_P$,
   $P_{\Delta'}=\cap \{P_i : v_i \in V_{\Delta'}\}\subset M'_1$
   and let  $C_i = \cup \{ C_{\Delta'}:\Delta' \subset \Delta_P,\dim \Delta' \geq i \}$.
    Note
  that $C_{\Delta'}$ is a singleton and $P_{\Delta'}=P$ for $\Delta'=\Delta_P$.
  Fix a point $x_P \in \Int P$ and let $W$ be an open neighborhood
   of $x$ in $M'_1$ such that $P_{\Delta'} \cap W$ is contractible
   for every $\Delta' \subset \Delta_P$.
  We will  construct an embedding
  $e_P : \Delta_P \lo M'_1$ such that
    $e_P(C_{\Delta'}) \subset P_{\Delta'}\cap W$ for every
   $\Delta' \subset \Delta_P$ and
   $e(C_{\Delta'} \setminus C_i) \subset
   \Int (P_{\Delta'} \cap W)$
   for every $i$ and  $\Delta'\subset \Delta_P$ such that $\dim \Delta' =i-1$.
   The map $e_P$ is constructed by induction on $i=t,t-1, \dots$ as follows.

   Set $e_P(C_{\Delta_P})=x_P \in \Int (P\cap W)$ and
  assume that $e_P$ is already defined on
  $C_i $
  such that the required properties are satisfied.
  Then
  for an  $(i-1)$-dimensional simplex  $\Delta' \subset \Delta_P$,
  $e_P(C_i \cap C_{\Delta'}) \subset \partial (P_{\Delta'}\cap W)$ and,
   since $P_{\Delta'}\cap W$ is contractible,
    $e_P$ can be extended over $C_{\Delta'}$
  such that
    $e_P(C_{\Delta'}) \subset P_{\Delta'}\cap W$ and
  $e_P(C_{\Delta'}\setminus C_i) \subset
 \Int (P_{\Delta'}\cap W)$. This way we extend $e_P$  over
  $C_{i-1}$ and the construction
  of $e_P$ is completed. Is is easy to see that $e_P$
  witnesses the intersection $P$. Since ${\cal P}_1$
  has no finite intersection of $\dim \leq l_1-2$,
   $\dim \Delta_P =t \leq m_1-l_1 +1=n$. Then replacing $e_P$
   by an approximation which also witnesses the intersection $P$
   we may assume that $e_P(\Delta_P) \subset M_1(k_1)$.
   Note that taking $W$ sufficiently close to $x_P$
   we can construct $e_P$ such  that the images of $e_P(\Delta_P)$  for
  all finite intersections of ${\cal P}_1$ will form
  a discrete family in $M'_1$.

  Fix a collection of maps $e_P$ with images contained in $M_1(k_1)$
   and forming a discrete family in $M'_1$
   that witnesses
  all the intersections of ${\cal P}_1$.    Extend $f : A_1 \lo A_2$
  to a closed embedding  $g_1 : X_1 \lo X_2$ such that $g_1(X_1)$
  is a $Z$-set in $X_2$.
  Approximate the map $ g{_1}^{-1}|_{\dots} : g_1(X_1) \lo X_1$
  by a closed
  embedding  $g_2 : g_1 (X_1) \lo X_1$ such that $g_2$ restricted to
  $A_2$ coincides with $f^{-1}$, $g_2(g_1(X_1))$ is a $Z$-set in
  $X_1$ and for every finite intersection $P$ of ${\cal P}_1$
  the map $g_2 \circ g_1 \circ e_P : \Delta_P \lo X_1 \setminus A_1 =M_1(k_1) \subset M'_1$
  also witnesses the
  intersection $P$.
   Finally extend $g_2$ to a
  closed embedding $g : X_2 \lo X_1$.

 Then one can find   an open subset  $M'_2$ of
   $M_2$ such that $M_2(k_2) \subset M'_2$ and $M'_2$ admits
   a rational triangulation ${\cal T}$ such that
    1) ${\cal T}$ properly approaches $A_2$,    2)
   for every intersection $P=P_0\cap  \dots \cap P_t$ of distinct
   elements $P_0,\dots, P_t$ of ${\cal P}_1$ and
   the map
   $g_1 \circ e_P : \Delta_P=[v_0,\dots, v_t] \lo M_2(k_2) \subset
   M'_2$ every simplex of $\cal T$ intersecting $g_1( e_P(\Delta_P))$
   does not intersect  $g^{-1}((M'_1 \setminus \Int (P_0 \cup \dots \cup P_t))(k_1))$,
   every simplex
   of ${\cal T}$ intersecting
   $g^{-1}(P_i(k_1))$ does not intersect $g_1(e_P(\Delta_P^i))$
   and
   every simplex of ${\cal T}$ containing $g_1(e_P(v_i))$ is
   contained in $g^{-1}((\Int P_i)(k_1))$ and
   3) for every finite intersections $P$ and $P'$ of ${\cal P}_1$
   such that $P \cap P' =\emptyset$, every simplex of $\cal T$
   intersecting $g^{-1}(P(k_1))$ does not intersect
   $g^{-1}(P'(k_1))$.

  Let $\cal B$ be the partition of $M'_2$
  formed
  by the stars of vertices of ${\cal T}$ with respect to the first
  barycentric subdivision of ${\cal T}$.
  Arrange the elements of  ${\cal P}_1 =\{ P'_1,P'_2,\dots \}$
  into a sequence and define $\mu(P'_1)=$ the union of the elements of
  ${\cal B}$ that intersect $g^{-1}(P'_1(k_1))$, $\mu(P_2)=$
  the union of the elements of ${\cal B}$ that intersect
 $g^{-1}(P'_2(k_1))$  but do not intersect $g^{-1}(P'_1(k_1))$,
 $\mu(P'_3)=$ the union of the elements of ${\cal B}$ that intersect
 $g^{-1}(P_3(k_1))$  but do not intersect $g^{-1}(P'_1(k_1))\cup g^{-1}(P'_2(k_1)) $
 and so on.  Denote ${\cal P}_2 =\{\mu(P'_1), \mu(P'_2),\dots \}$.
 By \ref{facts}, ${\cal P}_2$ is a partition of $M'_2$.

 Then the property 2)  guarantees that for every $P \in {\cal P}_1$,
 $g_1(e_P(\Delta_P)) \subset  \mu(P)$ and therefore $\mu(P) \neq
 \emptyset$  and for every non-empty intersection $P=P_0 \cap \dots \cap
 P_t$ of distinct elements $P_0,\dots P_t \in {\cal P}_1$, the
 map $g_1\circ e_P$ witnesses the intersection of
 $\mu(P_0), \dots ,\mu(P_t)$ and therefore
 $\mu(P_0) \cap \dots \cap \mu(P_t)\neq \emptyset$.
   The property 3) implies that we do not create
 additional intersections in ${\cal P}_2$ and hence
 $\mu : {\cal P}_1 \lo {\cal P}_2$ is a matching of
 partitions on $M'_1$ and $M'_2$ respectively. From the property 1)
 we conclude that ${\cal P}_2$ properly approaches $A_2$ and
 $\mu $ agrees with $f$.
 \\
 \\
 {\bf Improving connectivity of ${\cal P}_2$.}
 Without loss of generality we may
 replace $M_1$ and  $M_2$ by $M'_1$  and $M'_2$ respectively
 and  assume  that $M_1=M'_1$ and $M_2=M'_2$ for
 the  partition ${\cal P}_1$ and the initial partition ${\cal P}_2$.
 Let us show that using \ref{twin} we can modify ${\cal P}_2$
 into the required $l_2$-co-connected rational  partition. The construction
 of \ref{twin}  is  a combination  of \ref{total} and
 \ref{absorb}, and \ref{total} consists of \ref{simul}
 applied finitely many times.
 Abusing the notation we always  denote   by $M_2$, $F$ and ${\cal P}_2$ the manifold,
 the PL-subcomplex of $M_2$ and the decomposition
 of $M_2$
 which are the input of \ref{simul} and \ref{absorb}
 (=the output the previous applications of \ref{simul} and \ref{absorb})
 and  by $M'_2$, $F'$ and ${\cal P}'_2$ the output of
 \ref{simul} and \ref{absorb}
 (=the modifications   of $M_2$, $F$ and ${\cal P}_2$
  that are obtained
 after applying \ref{simul} and \ref{absorb}),
  and we denote again by $\mu$  the correspondence
  $\mu : {\cal P}_1 \lo {\cal P}'_2$  which is
  the composition of $\mu :{\cal P}_1 \lo {\cal P}_2$ with
   the natural
  one-to-one correspondence  between ${\cal P}_2$ and its
  modification
  ${\cal P}'_2$ (defined by sending each element
  of ${\cal P}_2$ to its modification).
 Thus we start with
 the initial partition partition ${\cal P}_2$ that during \ref{twin} turns into a decomposition
 forming a partition on $U=M_2 \setminus F$
 and,  gradually reducing the dimension of $F$,  ${\cal P}_2$
  returns to be  a  partition at the end of  \ref{twin}.
 Recall that after applying \ref{simul} and \ref{absorb}
 we always obtain an open subset $M'_2$ of $M_2$
 such that
 $R=M_2 \setminus M'_2$ is a PL-presented subset of $M_2$
 with $\dim R \leq n=m_2-l_2+1 \leq k_2$.
 By \ref{rational} we can choose
  a homeomorphism of $h: A_2 \cup M_2\lo A_2 \cup M_2$
  such that  $h$ is arbitrarily close to the identity map,
   $h(a)=a$ for $a \in A_2$, $h(R)$ is a rationally presented closed subset
  of $M_2$, $h(F')$ is a rational subcomplex of
  $h(M'_2)$,  $h({\cal P}'_2)$ is a rational decomposition of $h(M'_2)$
  and $h$ restricted to $M_2\setminus R$ is a PL-homeomorphism to
  $M_2\setminus h(R)$. Then
  replacing $M'_2$, $F'$  and ${\cal P}'_2$ by $h(M'_2)$, $h(F')$ and
  $h({\cal P}'_2)$ respectively we can always assume that
  $M_2(k_2) \subset M'_2$, $F'$ is a rational subcomplex of $M'_2$
   and ${\cal P}'_2$ is a rational decomposition.
  Thus at the end of \ref{twin} we get that $\dim F' \leq n=m_2-l_2+1 \leq k_2$
  and hence $M_2(k_1) \subset M'_2\setminus F'$. Now  we can replace
  the final modification $M'_2$ of $M_2$ by $M'_2 \setminus F$ and get
  the required  partition on $M'_2$.

  It is easy to see that  \ref{simul},  \ref{absorb} and \ref{rational}
  can be applied
 such that  each element
 of ${\cal P}_2$ will intersect its modification.
   This implies that $\mu : {\cal P}_1 \lo {\cal P}'_2$
   will automatically agree with $f$
  provided ${\cal P}'_2$ properly approaches $A_2$.
 Thus in order to get the partition required in  Proposition \ref{p1}
 the only thing we need to take care of is
  to perform  \ref{simul} and \ref{absorb} preserving
  the property that  the modification
   of  ${\cal P}_2$ properly approach $A_2$.

 Let us first analyze \ref{simul}.
   We want to carry out  \ref{simul}
 in such a way that
  the modification of ${\cal P}_2$ (the output  of \ref{simul}) will  properly
  approach $A_2$.
  Assume that ${\cal P}_2$ properly approaches $A_2$.
  Enlarge $\st({\cal P}_2,{\cal P}_2) $ to an open cover ${\cal C}$ of
   $M_2$ such that ${\cal C}$ properly approaches $A_2$
   and $\st({\cal P}_2,{\cal P}_2)$ refines ${\cal C}$.
   By Proposition  \ref{p2}
    there is an open cover ${\cal V}$
   of  $M_2(k_2)$ such that $\cal V$ properly approaches $A_2$,
    $\ {\cal C} | M_2(k_2)$
   refines $\cal V$ and for every $C \in {\cal C} $
   there is $V_C \in \cal V$ such that
   $C \cap M_2(k_2) \subset V_C$ and the inclusion
    $C\cap  M_2(k_2) \subset V_C$
    induces
  the zero-homomorphism
  of the homotopy groups in $\dim \leq n-1$.
  Then ${\cal W}=\st ( {\cal V}, {\cal C})$ is an open cover of $M_2$
  which properly approaches  $A_2$. Since every map from a space of $\dim \leq n$
  to an open subset $G$ of $M_2$ can be homotoped inside $G$
  to a map into $G\cap M_2(k_2)$ we conclude that
  the inclusion
    $C \subset \st (V_C, {\cal C}) \in \cal W$
    induces
  the zero-homomorphism
  of the homotopy groups in $\dim \leq n-1$ for every $C \in \cal C$.
  Recall that $\st({\cal P}_2,{\cal P}_2)$ refines ${\cal C}$. Then
   the construction \ref{simul}
  can be carried out  such that
  the output of \ref{simul} (the decomposition obtained
  after applying \ref{simul}) will refine  $\st({\cal W}, {\cal W})$
  and therefore will properly approach $A_2$.

  Now we will analyze \ref{absorb}.
  Once again we assume that ${\cal P}_2$
  properly approaches $A_2$ and we want to carry out  \ref{absorb}
 in such a way that
  the modification of ${\cal P}_2$ (the output  of \ref{absorb}) will  properly
  approach $A_2$ as well.
  The collection  $\st({\cal P}_2, {\cal P}_2)$ properly approaches
  $A_2$.  Since
  for $P \in {\cal P}_2$,  $\st(P, {\cal P}_2)$ is
  a union of elements of ${\cal P}_2$ we can  naturally define
  the set $\mu^{-1}(\st(P, {\cal P}_2))$ as the  union
  of the corresponding elements of ${\cal P}_1$.
   Then the collection
  $\mu^{-1}(\st({\cal P}_2, {\cal P}_2))$ properly
  approaches $A_1$
   since $\mu$ agrees with $f$.
  By proposition \ref{p2} there is an open cover
  ${\cal V}$ of $M_1(k_1)$ such that ${\cal V}$ properly
  approaches $A_1$ and for every $P \in {\cal P}_2$
  there is $V_P \in {\cal V}$ such that the inclusion
  $(\mu^{-1}(\st(P, {\cal P}_2)))(k_1) \subset V_P$ induces
  the zero-homomorphism of the homotopy groups in
  $\dim \leq n-1$.
  Denote $H_P=\st(V_P,\mu^{-1}(\st({\cal P}_2,{\cal P}_2)))$
  and ${\cal H} =\{ H_P : P \in {\cal P}_2 \}=
  \st({\cal V}, \mu^{-1}(\st({\cal P}_2,{\cal P}_2)))$.
  Clearly $V_P \subset H_P$, $\mu^{-1}(\st(P, {\cal P}_2)) \subset H_P$
  and ${\cal H}$ properly approaches $A_1$.
   Since $\mu^{-1}(\st(P, {\cal P}_2))$ is a rational submanifold
  of $M_1$
  of $\dim =m_1$,
  every map  from a space of $\dim \leq n$ into
   $\mu^{-1}(\st(P, {\cal P}_2))$ can be homotoped
   inside $\mu^{-1}(\st(P, {\cal P}_2))$  into
  $(\mu^{-1}(\st(P, {\cal P}_2))(k_1)$ and therefore
   the inclusion
  $\mu^{-1}(\st(P, {\cal P}_2)) \subset H_P$ induces
  the zero-homomorphism of the homotopy groups in
  $\dim \leq n-1$.
     Since each element $H_P$ of  ${\cal H}$
  is a union of elements of ${\cal P}_1$ we can define
  $\mu(H_P)$ and $\mu({\cal H})= \{\mu(H_P) : P \in {\cal P}_2 \}$.
  Note that $\mu({\cal H})$ properly approaches $A_2$ and
  $\st({\cal P}_2,{\cal P}_2)$ refines $\mu({\cal H})$.

   Assume that $F$ is a PL-subcomplex of $M_2$
  such that for $U=M_2 \setminus F$ the decomposition
  ${\cal P}_2$ forms on $U$  an $l$-co-connected partition,
  $l \geq l_2$ for which  $\mu : {\cal P}_1 \lo {\cal P}_2$
  becomes
  a matching when ${\cal P}_2$ is restricted to $U$. Use
  the notation of \ref{absorb} and let $\Delta$ be a simplex
  in $F$ that we need to absorb. Take an element $P$ of
  ${\cal P}_2$  such that $\Delta  \subset P$. Then by
  \ref{match} the inclusion $\st(P,{\cal P}_2)\cap U \subset
  \mu(H_P)\cap U$ induces the zero-homomorphism  of the
  homotopy groups in
  $\dim \leq m_2 -l$
  and therefore $S_\Delta$ can be contracted
  to a point inside $\mu(H_P)\cap U$ and
  even inside $\mu(H_P)\cap N$ because
  $\mu(H_P)\cap N$ is a deformation retract of  $\mu(H_P)\cap U$.
  Extend $\mu({\cal H})$ to an open cover $\cal W$ of $M_2$ such that
  $\cal W$ properly approaches $A_2$ and $\mu({\cal H})$ refines
  $\cal W$.
  Then \ref{absorb}
  can be carried out such that  the modification
  of ${\cal P}_2$  will refine $\st({\cal W}, {\cal W})$
  and therefore it will properly approach $A_2$.
  This completes the proof of the proposition.
    \hfill $\Box$

    \end{section}

    \begin{section}{Appendix}

    \subsection{Extending partitions in the black hole
    modification}
    \label{ext}
    We adopt  the notations and the assumptions  of \ref{black}.
      Recall that $T$ can be represented as
   the product  $T=B^{q+1} \times S$ where $S$ is identified with
  $a \times S$ for a  point $a$ in $\partial B^{q+1}$.
  We consider
   $B^{q+1}$ and $S$ as PL-embedded
  in Euclidean spaces. It induces the corresponding PL-embeding of
  $T$ in the product of  Euclidean spaces and
  we   refer to  these Euclidean spaces when properties of linearity are
  used, thus
  we say that a  simplex is linear if it is  a simplex (linearly spanned
  by its vertices) in the corresponding Euclidean space and
  a map of a linear simplex is linear if it is the linear extension
  of its values on the vertices.

   Denote by  $p_S: T \lo S$ and $p_B : T \lo B^{q+1}$  the projections.
   Let ${\cal T}_T$ be a triangulation of $T$ such that the simplexes
   of ${\cal T}_T$ are linear,
   $p_S$ and $p_B$ restricted to every simplex of ${\cal T}_T$  are linear
    and ${\cal T}_T$
     underlies $S$ and the simplexes of ${\cal T}$ contained in $T$.
   Let ${\cal T}'_T$ be a subdivision of ${\cal T}_T$ such that
    the simplexes of ${\cal T}'_T$
   are linear
   and for every
   simplex $\Delta$ of ${\cal T}_T$,
    $p_S(\Delta)$ is a subcomplex of $S$ with respect to ${\cal T}'_T$.

    Denote by ${\cal T}''_T$ the (first) barycentric subdivision of
    ${\cal T}'_T$ and let ${\cal B}_T$ be the partition of $T$
    formed by the stars of the vertices of ${\cal T}'_T$ with
    respect to ${\cal T}''_T$. Then ${\cal B}_S={\cal B}_T |S$
    is a partition of $S$.
    For a non-empty finite intersection $B$
    of ${\cal B}_S$ which is the  intersection of distinct
    elements  $B_0, \dots , B_t$ of ${\cal B}_S$  let
    the vertex $v_B$ of ${\cal T}''_T$   be the barycenter
    of the simplex of ${\cal T}'_T$ spanned by the vertices
    contained in $B_0 \cup \dots \cup B_t$.

    For a finite intersection $B$ of ${\cal B}_S$
    denote $K_B = p_S^{-1}(B)$, $N_B =p_S^{-1}(\partial B)$ and
    $M_B = K_B \cap \partial T$.
    Note that $K_B=B^{q+1} \times B$, $N_B=B^{q+1} \times \partial B$,
     $M_B=\partial B^{q+1} \times B$,
     $\partial K_B=N_B \cup M_B$ and
     $N_B \cap M_B=\partial B^{q+1} \cap \partial B$.

    Let us show that ${\cal P}_L$ restricted to
    $M_B\setminus F$   and ${\cal P}_L$
    restricted to $\partial M_B \setminus F$
    are
     partitions.
    Assume that  $B \in {\cal B}_S$ (that is $\dim B=\dim S$).
    Note that ${\cal P}_L$ restricted to $\Int M_B \setminus F_L$
    is a partition.
    Take a point
    $x\in \partial M_B$ and let
    $\Delta_x$ be the smallest  simplex of ${\cal T}_T$ containing
    $x$. Then $v_B \in p_S(\Delta_x)$ and let $y \in \Delta_x$ be
    such that $p_S(y)=v_B$. Let $G$ be a neighborhood of $x$
    in $\partial M_B$ such that $G$ is contained in the star of
    $y$ with respect to ${\cal T}_T$. Then $(z,t) \lo z(1-t) + ty$,
    $z \in G, 0\leq t < 1$ defines a PL-embedding of
    $ G \times [0,1)$ into $M_B$ such that $G \times [0,1)$
    is a neighborhood of $x$ in $M_B$ and for every
    simplex $\Delta \in {\cal T}_T$ we have that
    $\Delta \cap (G \times [0,1))= (G\cap \Delta)\times [0,1)$.
    Recall that ${\cal T}_T$ underlies ${\cal P}_L|\partial T$ and
     $F_L\cap \partial T$, and
    note that $G \times (0,1)$ is an open subset of
    $\partial  T$.  Therefore
    ${\cal P}_L$ restricted to $(G\times(0,1))\setminus F_L$ is a partition.
     All these facts together
    imply by \ref{facts} that ${\cal P}_L$ restricted
    $(G \times [0,1))\setminus F_L$ and
    ${\cal P}_L$ restricted  to $G\setminus F_L$
    are partitions and hence
    ${\cal P}_L$ restricted to
    $M_B\setminus F_L$  and  ${\cal P}_L$ restricted to $\partial M_B\setminus F_L$
    are
     partitions.

     Now assume that for a finite intersection
     $  B' \in {\cal B}_S$ we already showed that
     ${\cal P}_L$ restricted to
    $M_{ B'}\setminus F_L $  and  ${\cal P}_L$ restricted to
    $\partial M_{B'} \setminus F_L$
    are
     partitions and let $B$ be a finite intersection of
     ${\cal B}_S$ such that $B \subset \partial  B'$
     and $\dim B +1 =\dim  B'$.
      Then replacing  $S$ by  $\partial  B'$ and
      $T$ by $B^{q+1} \times \partial B'$ we can
      repeat  the above reasoning  to show that
      ${\cal P}_L$ restricted to
    $M_B \setminus F_L$  and  ${\cal P}_L$ restricted
    to $\partial M_B\setminus F_L$
    are
     partitions.
     Thus we have shown that for every finite
     intersection of ${\cal B}_S$,  ${\cal P}_L$ restricted to
    $M_B\setminus F_L$  and  to $\partial M_B \setminus F_L$
    are
     partitions.

      Let $B$ be a non-empty  finite intersection of ${\cal B}_S$.
     Denote by
     $V$ the set of vertices  $v$ of ${\cal T}''_T$
     such that $v \neq v_B$, $v \in \partial T$ and
      $[v,v_B]$ is a simplex of ${\cal T}''_T$.
     For each $v \in V$ choose a point $v'$ lying
     on the interval $[v,v_B]$ connecting $v$ and $v_B$
     as follows. If   $v \in S$  then set $v'=v$
   and note that in this case $v'\in \partial B$.
         If $v \notin S$ then choose $v'$ such that
     if  $v' \neq v_B, v' \neq v$ and $p_S([v', v]) \subset \Int B$.
     For a simplex  $\Delta $ of ${\cal T}''_T$
     spanned by vertices $v_B,v_1 \dots, v_k$, with
     $v_1,\dots v_k \in V$ denote by $\Delta'$
     the simplex spanned by the vertices $v_B, v'_1, \dots v'_k$
     and denote by $W_B$ the union of all such simplexes
     $\Delta'$.
      Then $W_B$ is a PL-ball
       which is a closed  neighborhood
     of $v_B$ in $M_B$ such that $W_B \cap S= M_B \cap S=B$
     and  for every $x \in \partial W_B$ we have that
     the interval $[v_B,x]$ is contained in $W_B$ and
     $[v_B , x] \cap \partial W_B =x$.
     Let $\pi^W_B : W_B \setminus v_B \lo \partial W_B$ be the radial projection.
     Using a PL-homeomorphism represent    $K_B$ as $ W_B \times [0,1]$
      such that
     $W_B$ is identified  with $W_B \times 0$ and extend $\pi^W_B$ to the
     radial projection
     $\pi_B : K_B\setminus v_B=(W_B \times [0,1])\setminus v_B
      \lo \partial K_B \setminus
     \Int W_B =(\partial W_B \times [0,1]) \cup (W_B \times 1 )$ ($\pi_B$
     is defined
     with respect to the linear structure of  $ W_B \times [0,1]$ induced
     by the linear structures of $W_B$ and $[0,1]$).

      Let us show that
      ${\cal P}_L$ restricted
       $M_B \setminus (F_L \cup \Int W_B) $  is a partition.
      Take $x \in \partial W_B \setminus S$. Then there are  an open neighborhood
      $G$  of $x$ in $\partial W_B \setminus S$ and
      $\epsilon >0$  such that
       $G \times [0, \epsilon )$ is embedded into $M_B\setminus (S \cup \Int W_B)$ by the map
      $(z,t) \lo z(1+t) - tv_B, z \in G$ such  that
      for every simplex $\Delta $ of ${\cal T}''_T$,
       $\Delta  \cap (G \times [0, \epsilon))=
      (\Delta \cap G)\times [0, \epsilon)$.
      Recall that ${\cal T}''_T$ underlies ${\cal P}_L$ and $F_L$ restricted
      to $\partial T$ and
       ${\cal P}_L$ restricted to
    $M_B\setminus F_L$  and  to $\partial M_B \setminus F_L$
    are
     partitions.
       Since $G \times (0,\epsilon)$
       is an open subset of $M_B$,
       ${\cal P}_L$ restricted to $(G \times (0,\epsilon))\setminus F_L$ is a partition and therefore
       by,
       \ref{facts}, ${\cal P}_L$ restricted to $G\setminus F_L$ is a partition.
       Then, again by \ref{facts},  ${\cal P}_L$
       restricted to $(G \times [0, \epsilon ))\setminus F_L$ is a partition and
       hence ${\cal P}_L$ restricted
       $M_B \setminus (F_L \cup \Int W_B) $  is a partition (recall that $S \subset F_L$).

    For every $A \in {\cal A}$
   and every  finite intersection
    $B$ of ${\cal B}_S$  define   sets
     $A^B$ and $A^{\partial B}$ as follows.
     For $B=\emptyset$ set
     $A^{\emptyset} =A$ if $A \cap T=\emptyset$ and
     $A^{\emptyset}=A\cup S$ if $ A \cap T\neq \emptyset$, and
     for $\dim B=0$ set  $A^{\partial B} =A^\emptyset$.
     Now by induction on $\dim B$ define $A^B$ and $A^{\partial B}$
     so that
     $A^B=A^{\partial B}\cup \pi^{-1}_B (A^{\partial B}
     \cap (\partial K_B \setminus \Int W_B))$ for $\dim B \geq 0$
     and $A^{\partial B}=\cup \{ A^{B'} : B' \subset \partial B \}$
     for $\dim B \geq 1$.
     Clearly $A^B$ and $A^{\partial B}$ are PL-subcomplexes of $M$.
     Note that
      if $A \cap T =\emptyset$ then $A^B=A^{\partial B}=A$
     for every $B$.

    From the construction of
      $W_B$ and $\pi_B$ it follows that for
        every simplex $\Delta \in {\cal T}'_T$  intersecting $ W_B$
        we have
        $\pi^{-1}_B (\Delta \cap \partial W_B)\setminus S=(\Delta \cap W_B)\setminus S$.
    Then it is easy to see that $S^{\partial B}=S^B=S$,
    $(A^{\partial B} \cap M_B) \setminus S=
    (A^B \cap M_B )\setminus S =(A\cap M_B)\setminus S$
    and hence
    $(A^{B_1} \cap A^{B_2}) \setminus S = A^{B_1 \cap B_2} \setminus S$,
     $({A_1^B} \cap  {A_2^B}) \setminus S = (A_1 \cap A_2)^{B} \setminus S$
         for every $A, A_1, A_2 \in {\cal A}$
     and  finite intersections $B,B_1, B_2$ of ${\cal B}_S$.

     Denote
     ${\cal P}^B =\{ P^B : P \in {\cal P}_L \}$ and
    ${\cal P}^{\partial B}=\{ P^{\partial B}: P \in {\cal P}_L \}$.
    Recall that $S \subset F_L$.
     Let us show by induction
     on $\dim B$ that
     ${\cal P}^B$ restricted to $K_B \setminus  F_L^B$ and
     ${\cal P}^{\partial B}$ restricted to
     $ N_B \setminus  F_L^{ \partial B}$
     are partitions.

      For $\dim B=0$, $N_B =\emptyset $
      and therefore ${\cal P}^{\partial B}$
      restricted to
       $ N_B \setminus  F_L^{ \partial B}$
       is vacuously a partition.
     Let $B$ be a finite intersection of ${\cal B}_S$
      such that ${\cal P}^{\partial B}$ restricted to
      $N_B \setminus  F_L^{\partial B}$
      is a partition.
      Denote $R_B =\partial K_B \setminus \Int W_B$.
    From the construction of $\pi_B$ it follows that
     $K_B \setminus  v_B $ is PL-homemorphic to  $R_B \times [0,1)$
      such that
      $R_B$ corresponds to $R_B \times 0$,
       ${\cal P}^B$ restricted
      to $K_B \setminus  v_B $ corresponds to
      the decomposition  $({\cal P}^{\partial B} |R_B)\times [0,1)$
      of $R_B \times [0,1)$ and
      $(F^B_L\cap K_B) \setminus  v_B $ corresponds to
      $(F_L^{\partial B} \cap R_B) \times [0,1)$.
     Note that $F_L^{\partial B} \cap M_B = F_L \cap M_B$ and
      ${\cal P}^{\partial B}$ restricted to $M_B \setminus S$ coincides
      with  ${\cal P}_L$ restricted $M_B \setminus S$ and recall
      that ${\cal P}_L$ restricted to $\partial M_B \setminus F_L$ and
      $ M_B \setminus ( F_L\cup \Int W_B)$
      are partitions.  Then, by \ref{facts},
      ${\cal P}^{\partial B}$ restricted to
     $ R_B \setminus  F_L^{ \partial B}$ is a partition
     and therefore  ${\cal P}^B$ restricted to $K_B \setminus  F_L^B$
     is also a partition.

      Now assume $B$ is a finite intersection of ${\cal B}_S$ such that
        $\dim B >0$ and for every finite intersection $B'$ of ${\cal B}_S$
        such that $B' \subset \partial B$
        we have that ${\cal P}^{B'}$ restricted to $K_{B'} \setminus  F_L^{B'}$
        is a partition. Then, by \ref{facts}, ${\cal P}^{\partial B}$ is a partition
        on $N_B  \setminus F_L^{\partial B}$ since
        $\{ K_{B'} : B'$ is a finite intersection of ${\cal B}_S$
        such that $B' \subset \partial B\}$ is a partition of $N_B$.
       The induction is completed.

        Thus we have shown that
        ${\cal P}^B$ restricted to $K_B \setminus  F_L^B$ is
        a partition for every finite intersection $B$ of ${\cal B}_S$.
        For $A \in \cal A$  denote  $A^\pi =\cup \{ A^B : B \in {\cal B}_S \}$
        and define ${\cal P}^\pi_L =\{ P^\pi : P \in {\cal P}_L \}$.
        Note that   $F^B_L =K_B \cap F^\pi_L$ and
        ${\cal P}^\pi_L$ restricted to $K_B$ coincides
        with ${\cal P}^B$.
        Then,
        by \ref{facts}, ${\cal P}^\pi_L$ is a partition on
        $T\setminus  F_L^\pi$ since
         $\{ K_{B} : B$ is a finite intersection of ${\cal B}_S \}$
         is a partition of $T$.
        It is obvious that ${\cal P}_L^\pi$ restricted to $L\setminus S$
         coincides  with ${\cal P}_L$ restricted to $L\setminus S$ and
         therefore ${\cal P}_L^\pi $ is a decomposition of $M$ which is a partition on
         $M\setminus  F_L^\pi$.

      Now let us show that for $A_1, A_2 \in \cal  A$,
      $(A_1 \setminus A_2) \setminus S$ is a strong deformation retract
      of
      $(A_1^\pi \setminus A_2^\pi)\setminus S$.
      To this end it is enough
      to show
      that for every finite intersection $B$ of ${\cal B}_S$,
   $(C'_1 \setminus C'_2)\setminus S$
    is a strong deformation retract
      of
      $(C_1 \setminus C_2)\setminus S$
      where $C'_i =A_i^{\partial B} \cap \partial K_B$
      and $C_i =A_i^B \cap K_B$. From the construction it follows
      that $(C_1 \setminus C_2)\setminus S$ can be topologically
      represented as
      $ C' \times [0,1)$ with
      $C'=(C'_1 \setminus C'_2) \setminus (S \cup \Int W_B)$
      such that   $(C'_1 \setminus C'_2) \setminus S$
      is identified with $(C' \times 0 ) \cup ((C' \cap \partial W_B) \times [0,1))$.
      Note that $C'$ is a space admitting a triangulation  for which  $ C' \cap \partial W_B$
      is a subcomplex of $C'$ and therefore
      $(C' \times 0 ) \cup ((C' \cap \partial W_B) \times [0,1))$ is
      a strong deformation retract of $C' \times [0,1)$.

      Let $\Delta \in {\cal T}_T$ be contained in $\partial T$.
      Since $p_S$ is linear on $\Delta$
      and $p_S(\Delta)$ is a subcomplex of ${\cal T}'_T$ we can conclude
      that $\dim \Delta \cap M_B \leq \dim \Delta +\dim B -\dim S$
      every finite intersection $B$ of ${\cal B}_S$.
      Then, since $\dim \Delta^B \setminus S
      \leq \dim (\Delta \cap M_B)\setminus S +1$
      for $\dim B=0$ and
      $\dim \Delta^B\setminus S  \leq \dim \Delta^{\partial B}\setminus S +1$
      for $\dim B >0$, we get that
       $\dim \Delta^\pi \setminus S \leq \dim \Delta \setminus S+1$.
      Thus $\dim (A^\pi \cap T)\setminus S  \leq \dim (A \cap \partial T)\setminus S +1$
      for every $A \in \cal A$.

      The remaining properties required in \ref{black} are easy to verify.

  \subsection{ Digging  holes for improving connectivity
    of an intersection}\label{hole1}
    Let $M$ be a triangulated $m$-dimensional manifold, $F$ a PL-subcomplex of $M$
    and $\cal P$ a decomposition of $M$ which is a partition on $M\setminus F$.
    Assume that $M$ is $(q-1)$-connected, $m \geq 2q+1$ and
    let $  S^{q-1} \subset   \Int M \setminus F$ be
   a  PL-embedded $q$-dimensional sphere  so that  $\cal P$ restricted to $S^{q-1}$
    is a partition and
     for every finite intersection $P'$ of $ \cal P$ with $P' \cap S^{q-1} \neq \emptyset$,
    $S^{q-1}\cap P' $ is
     properly embedded in $P' \setminus F$
     (that is $ S^{q-1} \cap \partial (P' \setminus F)  = \partial ( S^{q-1}\cap P' )$.
    Extend the embedding of $S^q$
    to a PL-embedding of $R^m=R^q \times R^{m-q} \subset \Int M$ such that
    $S^{q-1}=\partial B^q$ bounds a cube $B^q \subset R^q$ where $R^q$ and $R^{m-q}$ are identified with
    $R^{q}\times O \subset R^m$ and $O \times R^{m-q} \subset R^m$.

    Fix a triangulation $\cal T$ of $M$ which is the 2-nd barycentric subdivision
    of a tringulation that underlies $S^{q-1}$, $F$ and $\cal P$. By a simplicial
    neighborhood of   $A \subset M$ in a subcomplex $X$ of $M$ we mean the union of
    the simplixes of $\cal T$ intersecting $A$ and contained in $X$. We may assume
    that $\cal T$ is chosen so that
     the simplicial neighborhood of  $S^{q-1}$ in $M$ is contained in $R^m$.
    For every finite intersection $P'$ of $\cal P$ that meets $S^{q-1}$
    denote by  $G_{P'}$ the simplicial neighborhood of $P' \cap S^{q-1}$ in $P'$.
    Note that $G_{P''}=G_{P'} \cap P''$ if $P'' \subset P'$,  $G_{P'}$ is
     a regular neighborhood of $P'\cap S^{q-1}$ in $P' \setminus F$ and
     $P'\cap S^{q-1}$ is locally flat in $P' \setminus F$ because
     $\dim P' - \dim P'\cap S^{q-1} =m-(q-1)\geq q+2 $.

     For a block bundle $\xi$ over a cell complex $X$, written $\xi/X$,
       denote by
     $\sigma_\xi$ the cells of $X$, by $\beta_\xi $ the blocks of
     the total space $E(\xi)$ and by $(\sigma_\xi, \beta_\xi)$ the
     pairs such that $\beta_\xi$ is the block over $\sigma_\xi$,
     see  \cite{rs}.  We say that $\xi$ underlies
     $Y\subset E(\xi)$ if $Y$ is a union of blocks
     of $\xi$ and $\xi$ underlies a collection ${\cal Y}$ of
     subsets of $E(\xi)$ if $\xi$ underlies each $Y \in \cal Y$.

     For every finite intersection $P'$ of $\cal P$
     that meets $S^{q-1}$
     we are going to construct  by induction on $\dim P' \cap  S^{q-1}$
     a block bundle $\xi_{P'} $
     over $P' \cap S^{q-1}$ such that $E(\xi_{P'})=G_{P'}$ and
     ${\xi_{P'}}|P''\cap S^{q-1} =\xi_{P''}$ if $P'' \subset P'$.
     Obviously such a block bundle exists if $\dim P' \cap S^{q-1}=0$.
    Assume that for a finite intersection $P'$ of $\cal P$ a block bundle
    $\xi_{P''}$ is already constructed for every finite intersection
    $P''$ of ${\cal P}$ such that $P''\subset P'$ and $P''\neq P'$.
    Then $\xi_{P''}$ define the correponding block bundle $\xi_{\partial P'}$ over
    $\partial(P' \cap S^{q-1} )$ and
         by Theorem 4.3 of \cite{rs} this block bundle extends
         to a block bundle $\xi_{P'}$ over $P'$ such that $E(\xi_{P'})=G_{P'}$.

         Thus we have constructed block bundles $\xi_{P'}$ for every
         finite intersection $P'$ of $\cal P$ that meets $S^{q-1}$
         and these block bundles define the corresponding
         block bundle $\xi$ over $S^{q-1}$.

         Let $B_1^{m-q}$ be a cube in $R^{m-q}$. Fix  $\epsilon =1/2$ and for each
          cell $\sigma_\xi$ of $S^{q-1}$
         define the pair $(\sigma_\eta, \beta_\eta)$ with $\sigma_\eta=\sigma_\xi$
         and  $\beta_\eta = \gamma_\xi \times B_1^{m-q}$ where
           $\gamma_\xi= \{ x : x=ts,1-\epsilon \leq t\leq 1+\epsilon,
         s\in \sigma_\xi \} \subset R^q$. Then the pairs
              $(\sigma_\eta, \beta_\eta)$
         form a block bundle $\eta$ over $S^{q-1}$ such that
       $E(\eta)$
         underlies a regular neighborhood of $S^{q-1}$ in $R^m$.
         Hence by Theorem 4.4 of \cite{rs}
         there is a PL-homeomorphism of $R^m$
          realizing  an isomorphism of $\xi$ and $\eta$.
          Let $e : R^m \lo M$ be
          the composition of such
         PL-homeomorphism with the original embedding of  $R^m$ in $M$.
         Then replacing the original embedding  of  $R^m$ into $M$ by the embedding  $e$
            we may assume that   $\xi$ coincides with  $\eta$.

         Fix a triangulation ${\cal T}_R$
         of $R^m$  such that the simplexes of  ${\cal T}_R$  are linear
         in $R^m$ and  ${\cal T}_R$  underlies  ${\cal P}|R^m$,
          $F\cap R^m$ and  $ B_1^q \times B_1^{m-q}$
          where  $B_1^q=(1+\epsilon)B^q$.
          Let $p : R^m =R^q \times R^{m-q} \lo R^{m-q}$ be the projection and let
          $a \in \Int B_1^{m-q}$ be such that
        $a \notin p(\Delta)$  for every $\Delta \in {\cal T}_R$
        with $\dim p(\Delta)< m-q$ and
        $a \notin \partial p(\Delta)$  for every $\Delta \in {\cal T}_R$
        with $\dim p(\Delta)= m-q$. Take a cube $B_2^{m-q}$ such that
        $a+B_2^{m-q} \subset \Int B_1^{m-q}$ and
        $a+B_2^{m-q} \subset \Int p(\Delta)$ for every
        $\Delta \in {\cal T}_R$ such that
         $\Delta \cap ( B_1^q \times a) \neq \emptyset$.
         Then  the pairs $(\sigma_\tau, \beta_\tau)$ with
          $\sigma_\tau = (B_1^q \times a)\cap \Delta$
          and
         $\beta_\tau = (B_1^q \times (a+ B_2^{m-q}))\cap \Delta$,
          $\Delta \in {\cal T}_R$
         form a block bundle $\tau$ over $B_1^q \times a$ and  we assume that
         $B_2^{m-q}$ is so small that
         $E(\tau)= B_1^q \times (a+ B_2^{m-q})$ and
         $E(\tau | \partial B_1^q) =\partial B_1^q \times (a +B_2^{m-q})$.

         By Theorem 1.1 of \cite{rs} every block bundle over a cube is trivial.
         Hence $\tau$ is isomorphic to
         a trivial block bundle  $\delta$ over $B_1^q$,
         that is there is a PL-homeomorphism $h: E(\delta) \lo E(\tau)$
         where $E(\delta)$ is the product
          $E(\delta)= B_1^q \times B_3^{m-q} \subset R^m$ of $B_1^q$
          with
          a cube $B_3^{m-q}$ in $R^{m-q}$ such that
           $h(x)=(x,a)$ for every $x \in B_1^q$,
          $(\sigma_\delta, \beta_\delta)$ is a pair of $\delta$
          if and only if
            $(h(\sigma_\delta), h(\beta_\delta))$ is a pair of $\tau$
            and, for every pair
          $(\sigma_\delta, \beta_\delta)$  of $\delta$,
             $\beta_\delta = \sigma_\delta \times B_3^{m-q}$.

          Identify $B_1^q \times B_3^{m-q}$ with
          $h(B_1^q \times B_3^{m-q})$ and
            let $B^{m-q}$ be a cube lying in  $\Int B_3^{m-q}$.
            Then from the construction it follows that

            (*)
            $P' \cap (B_1^q \times \Int B_3^{m-q}) =
           (P' \cap B_1^q) \times \Int B_3^{m-q}$
           for every
             $P' \in {\cal P}$ and
        $F \cap (B_1^q \times \Int B_3^{m-q}) =
           (F \cap B_1^q) \times \Int B_3^{m-q}$.

           From the construction it also follows that
           $\cal P$ restricted to the sets $\Int B_1^q \setminus \Int B^q$,
           $\partial B^q$ and $B^q \setminus F$ are partitions.
           Hence by (*) and  \ref{facts}, $\cal P$ restricted to
           the sets
           $(\Int B_1^q \setminus \Int B^q) \times \Int B_3^{m-q}$,
           $(B^q \setminus F) \times (\Int B_3^{m-q} \setminus \Int B^{m-q})$,
            $\partial B^q \times (\Int B_3^{m-q} \setminus \Int  B^{m-q})$,
            $\partial B^q  \times \partial B^{m-q}$,
            $\partial B^q \times B^{m-q}$ and
            $(B^q \setminus F)\times \partial B^{m-q}$ are partitions.
            Then, once again by \ref{facts}, $\cal P$ restricted
            to $\Int (B_1^q \times B_3^{m-q}) \setminus (F \cup  (\Int(B^q \times B^{m-q})))$
            and $\partial (B^q \times B^{m-q})\setminus F$ are partitions.

            Denote $B^m=B^q \times B^{m-q}$. Thus we get that

            (**)  $\cal P$ restricted to $M \setminus (F \cup \Int B^m)$
            and $\partial B^m\setminus F$ are partitions.

            Let $P'$ be a finite intersection of $\cal P$.
            Let us show
            that

            (***) $ P'\setminus (F \cup \Int B^m)$ is
            a  deformation retract of  $ P'\setminus (F\cup \Int B^q)$.

            Note the following general fact: for a triangulated space $A_1$,
            a subcomplex $A_2$ of $A_1$ and a cube $B$  the space
           $(A_2 \times B)\cup (A_1 \times \partial B)$ is a strong
            deformation retract of  $(A_2 \times B)\cup (A_1 \times ( B\setminus O))$.
            Then for $A_1=(P'\cap B^q)\setminus F$ and
            $A_2=P'\cap \partial B^q$ we have by (*) that
            $B^m \cap (P'\setminus (F \cup \Int B^m))=
            (A_2 \times B^{m-q})\cup (A_1 \times \partial B^{m-q})$
           and
          $B^m\cap (P' \setminus (F \cup \Int B^q))=
                 (A_2 \times B^{m-q})\cup (A_1 \times ( B^{m-q}\setminus O))$.
       Thus  $B^m \cap (P'\setminus (F \cup \Int B^m))$ is a strong deformation
        retract of $B^m\cap (P' \setminus (F \cup \Int B^q))$ and hence
        (***) holds.

                 Now let us compare the embedding of $S^{q-1}= \partial B^q$
           into $M$ under $h$ with the embedding of $S^{q-1}$ under  $e$
           (which coincides with the original embedding of $S^{q-1}$
           into $M$). One can easily observe that
            there is a PL-homeomorphism
           $\psi : M \lo M$ such that
           $\psi\circ h|S^{q-1}=e|S^{q-1}$,
              $\psi$ does not move the points outside $E(\eta)$,
            $\psi(P')=P'$ and
           for every $P' \in \cal P$. Since
            $\psi$ does not move the points outside $E(\eta)$
            we have that
            $\psi(F)=F$.  Then
            replacing the embedding of $B_1^q \times B_3^{m-q}$
           into $M$ by its composition with $\psi$  we preserve the properties
            (*), (**) and (***) and obtain that the embedding of $S^{q-1}$ coincides
           with the original embedding of $S^{q-1}$.

           Thus we can extend  the original embedding of $S^{q-1}$
           to
           a PL-embedding  of $B^m =B^q\times B^{m-q}$ such that
           $\cal P$ restricted to $\partial B^m \setminus F$ and
           $M \setminus (\Int B^m  \cup F)$ are
            partitions, $S^{q-1}= \partial B^q$,
            $F\cap B^m= (F\cap B^q)\times B^{m-q}$,
            $P'\cap B^m= (P'\cap B^q)\times B^{m-q}$ for every $P' \in \cal P$
            and $ P'\setminus (F \cup \Int B^m)$ is
            a  deformation retract of  $ P'\setminus (F\cup \Int B^q)$
            for every finite intersection $P'$ of $\cal P$.
            It is clear that $B^m \subset R^m \subset \Int M$.

            Note that if $S^{q-1}$ can be contracted to a point in
            an open subset $W$ of $M$ then we may assume that
            $B^m$ is contained in $W$.

   \subsection{Digging holes for absorbing simplexes}
   \label{hole2}
   We adopt the notations and assumptions of \ref{absorb}.
     Recall that $\cal T$ is the 2-nd barycentric subdivision of some triangulation
      ${\cal T}_0$ of $M$  underlying $\cal P$ and
     consider $M$ as embedded in a Hilbert space by an embedding
     which is linear on every simplex of
      ${\cal T}_0$.

   Fix $\Delta \in {\cal T}_F$ and let  $S^{q-1}=S_\Delta$ be the link of
   $\Delta$ with respect to $\cal T$,  $S^{m-1} =\partial (\Delta * S_\Delta)$
   and $b$ the barycenter of $\Delta$.

   Let us show that ${\cal P}$ restricted to $N$, $S^{m-1}\setminus \partial \Delta$
   and $(\Delta * S_\Delta)\setminus  \Delta$ are partitions.
   For a point $x \in S^{m-1}\setminus \partial \Delta$ choose an $\epsilon >0$
   and a neighborhood $G$ of $x$ in $S^{m-1}\setminus \partial \Delta$ such that
   the map $(z,t)\lo z+tb$, $z \in G$, $t \in (-\epsilon, \epsilon)$ embeds
   $G \times (-\epsilon, \epsilon)$ into $M \setminus F$
   and for every simplex $\Delta'$ of ${\cal T}_0$ that intersects
   $G \times (-\epsilon, \epsilon)$ we have that
   $\Delta' \cap (G \times (-\epsilon, \epsilon))=
   (\Delta' \cap G) \times (-\epsilon, \epsilon) $.
   Then for every $P \in \cal P$ that intersects
   $G \times (-\epsilon, \epsilon)$ we have that
   $P \cap (G \times (-\epsilon, \epsilon))=
   (P \cap G) \times (-\epsilon, \epsilon) $.
   Therefore, by \ref{facts},
   ${\cal P}$ restricted to $G$, $G \times (-\epsilon,0]$ and
   $G \times [0, \epsilon)$ are partitions and hence
   ${\cal P}$ restricted to $S^{m-1}\setminus \partial \Delta$,
   $(\Delta * S_\Delta)\setminus  \Delta$ and $N$ are partitions as well.

   Denote by  ${\cal T}_\Delta$ the collection of the simplexes of ${\cal T}$
   intersecting  $S^{q-1}$ and containing $\Delta$.
   Let ${\cal T}'$ be the second
   barycentric subdivision of ${\cal T}$. For every $\Delta' \in {\cal T}_\Delta$
   denote by $G_{\Delta'}$ the simplicial
   neighborhood of $\Delta'  \cap S^{q-1}$ in $\Delta' \cap S^{m-1}$
   with respect to ${\cal T}'$.
   In a way similar to constructing $\xi_{P'}$ in \ref{hole1} we construct for every
   $\Delta' \in {\cal T}_\Delta$ a block bundle $\xi_{\Delta'}$ over
   $\Delta'  \cap S^{q-1}$ such that $E(\xi_{\Delta'})=G_{\Delta'}$ and
   $\xi_{\Delta'}|({\Delta''}\cap S^{m-1})=\xi_{\Delta''} $ if ${\Delta''} \subset \Delta'$.
   Then the block bundles $\xi_{\Delta'}$  define the corresponding
   block bundle $\xi$ over $S^{q-1}$ with $E(\xi)$ being the simplicial
   neighborhood of $S^{q-1}$ in  $S^{m-1}$ with respect to ${\cal T}'$.
   Note that $E(\xi)$ is a regular neighborhood of $S^{q-1}$ in $\partial N$.

   Extend the embedding of $S^{q-1}$ in $N$ to  a  PL-embedding
   of a $q$-dimensional cube $B_1^q$ into
   $N$ such that $S^{q-1}= \partial B_1^q$ and $B_1^q \cap \partial N =S^{q-1}$
   and note that the embedding of $B_1^q$ is locally flat because $m \geq 2q+1$.
   Extend the regular neighborhood $E(\xi)$ of $S^{q-1}$ in $\partial N$ to
   a regular neighborhood of $B_1^q$ in $N$ and by Theorem 4.3 of \cite{rs}
   represent this neighborhood as the total space $E(\nu)$ of a block  bundle
   $\nu$ over $B_1^q$  such that $\nu | S^{q-1} = \xi$. Since $\nu$ is
   a trivial block bundle, $E(\nu)$ is PL-homeomorphic to
   the product
   $E(\nu)=B_1^q \times B_1^{m-q}$ of $B_1^q$ with an $(m-q)$-dimensional
   cube $B_1^{m-q}$ such that for each pair
   $(\sigma_\nu, \beta_\nu)$ of $\nu$, $\beta_\nu= \sigma_\nu \times B_1^{m-q}$.
   Fix a triangulation ${\cal T}_\nu$ of $E(\nu)$ which underlies ${\cal T}$
   restricted to $E(\nu)$.
   In a way similar to constructing the block bundle $\tau$ in \ref{hole1}
   choose
   a point $a$  in $B_1^{m-q}$  and a cube $B_2^{m-q}$
   such that $a +B_2^{m-q}\subset \Int B_1^{m-q}$
   and the pairs
    $(\sigma_\tau, \beta_\tau)$ defined by
     $\sigma_\tau = \Delta' \cap (B_1^q \times a)$,
     $\beta_\tau = \Delta' \cap (B_1^q \times (a +B_2^{m-q}))$, $\Delta' \in {\cal T}_\nu$
     form a block bundle $\tau$ over $B_1^q \times a$ with
     $E(\tau)= B_1^q \times (a +B_2^{m-q})$ and
     $E(\tau | \partial (B_1^q \times a))=
     E(\tau) \cap (\partial B_1^{q} \times (a+ B_2^{m-q}))
     (=E(\tau)\cap S^{m-1}$).
     Then the block bundle $\tau$ underlies $\cal P$ and  ${\cal T}_\Delta$ restricted
     to $E(\tau)$.

     Note that since the block bundle $\nu$ is trivial,
     the block bundle $\xi$ is trivial as well and
     $E(\xi)=\partial B_1^q \times B_1^{m-q}$. Then
     it is obvious that
     the projection  $\partial B_1^q \times a \lo \partial  B_1^q$ extends to
     a PL-homeomorphism
     $h : E(\xi) \lo E(\xi)$ such that $h(x)=x $ for every $x \in  \partial E(\xi)$,
      $h(\beta_\xi)=\beta_\xi$ for every  block $\beta_\xi$ of $\xi$ and there is
      a PL-isotopy $H :  E(\xi)\times [0,1]\lo E(\xi)\times [0,1]$ relative to
      $\partial  E(\xi)$      so that $H_0=h$.
     Recall that    the barycenter of
     $\Delta$ is denoted by $b$. Then, since $\partial \Delta  \cap E(\xi) =\emptyset$,
     there is $\epsilon >0$ such that
       for
     every $x \in E(\xi)$, $P \in \cal P$ and $t\in [0, 1]$
     the point $b+(1+t\epsilon)(x-b)$ belongs to $P$ if and only if $x$ belongs to $P$.
      Embed $E(\xi)\times [0,1]$ into $N$ by sending $(x,t)$ to $b+(1+t\epsilon)(x-b)$
      and define the PL-homeomorphism  $g : N \lo N$ such that $g$
      coincides with $H$ on $E(\xi)\times [0,1]$ and $g$ does not move the points
      outside $E(\xi)\times [0,1]$. From the properties described above it is
      clear that $g(\partial (B_1^q \times a))=\partial B_1^q$
      and  the block bundle  $g(\tau)$ underlies $\cal P$ and ${\cal T}_\Delta$
      restricted to $E(g(\tau))$.
      Denote $B_2^q =g(B_1^q \times a)$ and $\theta=g(\tau)$.

      Thus we get
      a block bundle $\theta$ over a cube $B_2^q$ such that
      $S^{q-1}=\partial B_2^q= \partial N  \cap B_2^q$,
      $\theta$ underlies $\cal P$ and ${\cal T}_\Delta$ restricted to $E(\theta)$
      and $E(\theta|\partial B_2^q)= S^{m-1}\cap E(\theta)$.
      Then, since $\theta$ is trivial, $E(\theta)$ is PL-homemorphic to
      $B_2^q \times B_3^{m-q}$ such that
            $\partial N \cap (B_2^q \times B_3^{m-q})=
            \partial B_2^q \times B_3^{m-q} \subset S^{m-1}$
            and
      for  each pair
      $(\sigma_\theta, \beta_\theta)$ of $\theta$,
      $\beta_\theta= \sigma_\theta \times B_3^{m-q}$.
      Take a cube $B_4^{m-q}$ such that $ B_4^{m-q} \subset
      \Int B_3^{m-q}$ and define the block bundle $\rho$
      over $B_2^q$
      by  restricting the blocks of  $\theta$ to
      $E(\rho)=B_2^q \times B_4^{m-q}$. Using the reasoning similar
      to the one applied for proving the property (**) in
      \ref{hole1} we can show that ${\cal P}$ restricted
      $N \setminus \ (B_2^q \times \Int B_4^{m-q})$ and
      $\partial ((\Delta* S_\Delta) \cup (B_2^q \times  B_4^{m-q}))\setminus \partial \Delta$
      are partitions. It is clear that $E(\rho) \subset \Int M$.

      Take a pair  $(\sigma_\rho^\Delta, \beta_\rho^\Delta)$ of $\rho$
      such that $\dim \beta_\rho^\Delta=m$
      and let $P_\Delta \in\cal P$ be such that
      $\beta_\rho^\Delta \subset P_\Delta$. Since $\rho$ is a trivial block bundle
     we can replace $\rho$
      by another block bundle subdividing the cells of $B_2^q$
      into smaller cells and defining the corresponding blocks using the product
      structure. This way we may assume that cells of $B_2^q$ are so small that
      the pair $(\sigma_\rho^\Delta, \beta_\rho^\Delta)$  can be chosen so that
          $\beta_\rho^\Delta \subset \Int P_\Delta$.

      Let the cube $B_3^q \subset \Delta * S_\Delta $
      be the join $B_3^q=b*S_\Delta$.
      Define a block bundle $\psi$ over $B_3^q$ as follows:
      the pairs  $(\sigma_\rho, \beta_\rho)$  of $\rho$ with
      $\sigma_\rho \subset S^{q-1}$
       are    pairs of $\psi$,
       $(b, \Delta)$ is a pair of $\psi$
         and the rest of the pairs $(\sigma_\psi, \beta_\psi)$ of
       $\psi$ are of the form $\sigma_\psi=\Delta' \cap B_3^q$,
       $\beta_\psi=\Delta'$, $\Delta' \in {\cal T}_\Delta$.
       One can easily check that $\psi$ is indeed a block bundle
       and $E(\psi)=\Delta*S_\Delta$.
       Since $\psi$ and $\rho$ over $S^{q-1}$ coincide, they
       form the corresponding block  bundle $\phi$ over the sphere
       $B_3^q \cup B_2^q$. Note that $B^q=(B_2^q \cup B_3^q)\setminus \Int \sigma^\Delta_\rho$
       is a $q$-dimensional cube. Define the block bundle $\eta$ as $\psi$ restricted to
       $B^q$. Then
       $(b,\Delta)$ is a pair of $\eta$, $b \in \Int B^q$,
       $\partial B^q = \partial \sigma_\rho^\Delta $,
       $E(\eta | \partial B^q) \subset \Int P_\Delta$,
        $\eta$ underlies $\cal P$ restricted to $E(\eta)$,  $E(\eta)\cap F=\Delta$ and hence
        $\eta$ underlies $F$ restricted to $E(\eta)$.

        Thus assuming that the center $O$  of $B^q$ is located at $b$ we can represent
        $E(\eta)$ as the product $B^m=B^q \times B^{m-q}$ such that
        $\Delta=0\times B^{m-q}=F \cap B^m$,
        $\partial B^q \times B^{m-q} \subset \Int P_\Delta$ and
        $P \cap B^m=(P\cap B^q)\times B^{m-q}$ for every $P \in \cal P$ that
        meets $B^m$.
        By the reasoning similar to the one applied to the proof of (***)  in \ref{hole1}
        these properties imply that
         for every finite intersection
        $P$ of $  \cal P$,
        $ P\setminus (F \cup \Int B^m)$ is
            a  deformation retract of  $ P\setminus (F\cup \Int B^q)$.
           Since  $B^m=(\Delta*S_\Delta)\cup ((B_2^q \setminus
            \Int \sigma^\Delta_\rho )\times B_4^{m-q})$
            and $\beta_\rho^\Delta=\sigma^\Delta_\rho \times B_4^{m-q}$ is
            contained in $\Int P_\Delta$ we derive from the properties
            of $\rho$  that
       $\cal P$ restricted to $M\setminus(F\cup \Int B^m)$
        and $\partial B^m \setminus F$ are partitions and $B^m \subset \Int M$.

          Note that if $S^{q-1}$ can be contracted to a point in
            an open subset $W$ of $N$ then the construction can be carried out such that
           $B^m$ is contained in $(\Delta * S_\Delta) \cup  W$.

    \end{section}

Department of Mathematics\\
Ben Gurion University of the Negev\\
P.O.B. 653\\
Be'er Sheva 84105, ISRAEL  \\
e-mail: mlevine@math.bgu.ac.il\\\\
\end{document}